\documentclass[article]{elsarticle}

\usepackage{lineno,hyperref}
\usepackage{color}
\usepackage{amsfonts}
\usepackage{amsmath}
\usepackage{amssymb}
\usepackage{bm}
\usepackage{mathabx} 
\usepackage{epsfig}
\usepackage{multirow}
\usepackage{subcaption}
\usepackage{graphicx}
\usepackage{rotating}
\usepackage{listings}
\usepackage{tikz, pgf}
\usepackage{tikz-3dplot}
\usepackage[normalem]{ulem}

\makeatletter
\def\@xfootnote[#1]{%
\protected@xdef\@thefnmark{#1}%
\@footnotemark\@footnotetext}
\makeatother

\newcommand{\bff}{\mathbf{f}}

\newcommand{\bfk}{\mathbf{k}}

\newcommand{\bfx}{\mathbf{x}}

\newcommand{\bfK}{\mathbf{K}}

\newcommand{\mm}[1]{\rm mm}

\newcommand{\beq}{\begin{equation}}
\newcommand{\eeq}{\end{equation}}
\newcommand{\bea}{\begin{eqnarray}}
\newcommand{\eea}{\end{eqnarray}}

\newcommand{\bit}{\begin{itemize}}
\newcommand{\eit}{\end{itemize}}
\newcommand{\ben}{\begin{enumerate}}
\newcommand{\een}{\end{enumerate}}

\newcommand{\bK}{\mathbf{K}}

\newcommand{\bx}{\mathbf{x}}

\modulolinenumbers[5]

\journal{Communications in Applied Mathematics and Computational Science}

\begin{document}
\begin{frontmatter}

\title{An Application of Gaussian Process Modeling for 
High-order Accurate Adaptive Mesh Refinement Prolongation}

\author[a]{Steven I. Reeves}
\ead{sireeves@ucsc.edu}
\author[a]{Dongwook Lee\corref{mycorrespondingauthor}}
\ead{dlee79@ucsc.edu}
\author[b]{Adam Reyes}
\ead{acreyes@uchicago.edu}
\author[c]{Carlo Graziani}
\ead{carlo@mcs.anl.gov}
\author[b,d,e]{Petros Tzeferacos}
\ead{p.tzeferacos@rochester.edu}

\cortext[mycorrespondingauthor]{Corresponding author}


\address[a]{Department of Applied Mathematics, 
The University of California, Santa Cruz, CA, United States}
\address[b]{Flash Center for Computational Science, Department of Astronomy \& Astrophysics, The University of Chicago, Chicago, IL, United States}
\address[c]{Mathematics and Computer Science, 
Argonne National Laboratory, Argonne, IL, United States}
\address[d]{Department of Physics and Astronomy, University of Rochester, Rochester, NY, United States}
\address[e]{Laboratory for Laser Energetics, University of Rochester, Rochester, NY, United States}
\begin{abstract}
We present a new polynomial-free 
prolongation scheme for Adaptive Mesh Refinement (AMR) 
simulations of compressible and incompressible computational fluid dynamics. 
The new method is constructed using a
multi-dimensional kernel-based Gaussian Process (GP) prolongation model.
The formulation for this scheme
was inspired by the GP methods introduced by A. Reyes \textit{et al.}
[A New Class of High-Order Methods for Fluid Dynamics Simulation 
using Gaussian Process Modeling, Journal of Scientific Computing, 76 (2017), 443-480; 
A variable high-order shock-capturing finite difference method with GP-WENO, 
Journal of Computational Physics, 381 (2019), 189–217].
In this paper, we extend the previous GP interpolations/reconstructions to a new GP-based AMR prolongation 
method that delivers a high-order accurate prolongation of data from
coarse to fine grids on AMR grid hierarchies.
In compressible flow simulations 
special care is necessary to handle shocks and discontinuities in a stable manner.
To meet this, we utilize 
the shock handling strategy using the GP-based smoothness indicators
developed in the previous GP work by  A. Reyes \textit{et al.}.
We demonstrate the efficacy of the GP-AMR method 
in a series of testsuite problems using 
the AMReX library, in which the GP-AMR method has been implemented. 
\end{abstract}

\begin{keyword}
Adaptive Mesh Refinement; 
Prolongations;
High-order methods;
Gaussian processes; 
Computational fluid dynamics;
\end{keyword}

\end{frontmatter}

\nolinenumbers


\section{Introduction}
\label{sec:introduction}

In the fields of geophysics, astrophysics, and laboratory plasma astrophysics simulations have become essential to characterizing and understanding complex processes
(e.g., ~\citep{glatzmaiers1995three,
jordan2008three,tzeferacos2015flash,meinecke2014turbulent}). 
As increasingly more complex systems to be considered for better computer modeling,
modern simulation codes face increasingly versatile challenges to meet expected metrics in a 
possibly vast parameter space.
CFD has been (and will continue to be) an indispensable tool to improve
our capabilities to investigate conditions where simplified theoretical models inadequately capture the correct physical behavior and experiments can be prohibitively expensive or too observationally difficult to be the sole pathways for discovery.
In these simulations flow conditions can develop in which the physics becomes extremely challenging to simulate
due to significant imbalances in length and temporal scales. 
To alleviate such conditions in computer simulations, practitioners have explored
approaches by which a computer simulation can focus on localized
flow regions when the dynamics exhibit confined features that evolve
on a much shorter length scale relative to the flow dynamics on the rest of the 
computational domain. 

Adaptive mesh refinement (AMR) is one such approach
that allows a local and dynamic change in the grid resolutions of a simulation
in space and time.
Since the 1980s, AMR has been an exceptional tool 
and has become a powerful strategy in utilizing computational fluid dynamics (CFD) simulations
for computational science across many disciplines 
such as astrophysics, geophysics, atmospheric sciences, oceanography, biophysics,
engineering, and many others \cite{plewa2005adaptive}.

There have been many advancements in AMR since the seminal paper by Berger and Oliger~\cite{amr_orig}.
In their paper, the primary concern was to focus on a strategy for generating subgrids and managing the grid hierarchy
for scalar hyperbolic PDEs in one and two spatial dimensions (1D and 2D).
In the subsequent work by Berger and Colella \cite{Berger}, further improvements were made possible
for numerical solutions of the 2D Euler equations to provide a robust shock-capturing AMR algorithm
that satisfies the underlying conservation property 
on large-scale computer architectures. The novel innovations in their work have now become the AMR standards,
namely including refluxing (or flux correction) between fine-coarse interface boundaries, 
conservative (linear) prolongation and restriction on AMR hierarchies, and
timestep subcycling.
Bell \textit{et al.} extended the precedent 2D AMR algorithms of~\cite{amr_orig,Berger} to a 3D
AMR algorithm and applied it to solve 3D hyperbolic systems of conservation laws~\cite{amr3d}.
They demonstrated that the AMR algorithm reduced the computational cost
by more than a factor of 20 than on the equivalent uniform grid simulations in simulating a 3D
dense cloud problem interacting a Mach 1.25 flow on Cray-2.
This is, by far, the main benefit of using AMR, particularly in large 3D simulations, 
in that one could gain such a computational speed-up 
by focusing computational resources on the dynamically interesting regions of the simulation.

AMR can be expected to become more computationally expensive relative to a uniform grid solution using high-order (4th or higher) PDE solvers, as a significant fraction of the computational domain becomes dominated by small scale structures .
Jameson \cite{jameson2003amr} estimated that small scale features, such as shocks or vortices,  should not exceed more than a third of the computational domain in order for low order AMR schemes to be computationally competitive. The effectiveness of AMR in atmospheric simulations has also been studied \cite{ferguson2016analyzing}.

Modern AMR implementations may be categorized into two main types: structured and unstructured. 
Unstructured AMR, and meshes in general, are very useful for problems with irregular geometry 
(e.g., many structural engineering problems), 
but is often computationally complex and difficult to handle when regridding. 
On the other hand, structured AMR (SAMR, or block-structured AMR)
offers practical benefits (over unstructured) such as ease of discretization, a global index space, 
accuracy gain through cancellation terms, and ease of parallelization. 

\begin{figure}[ht!]
\begin{center}
\includegraphics[width=11cm]{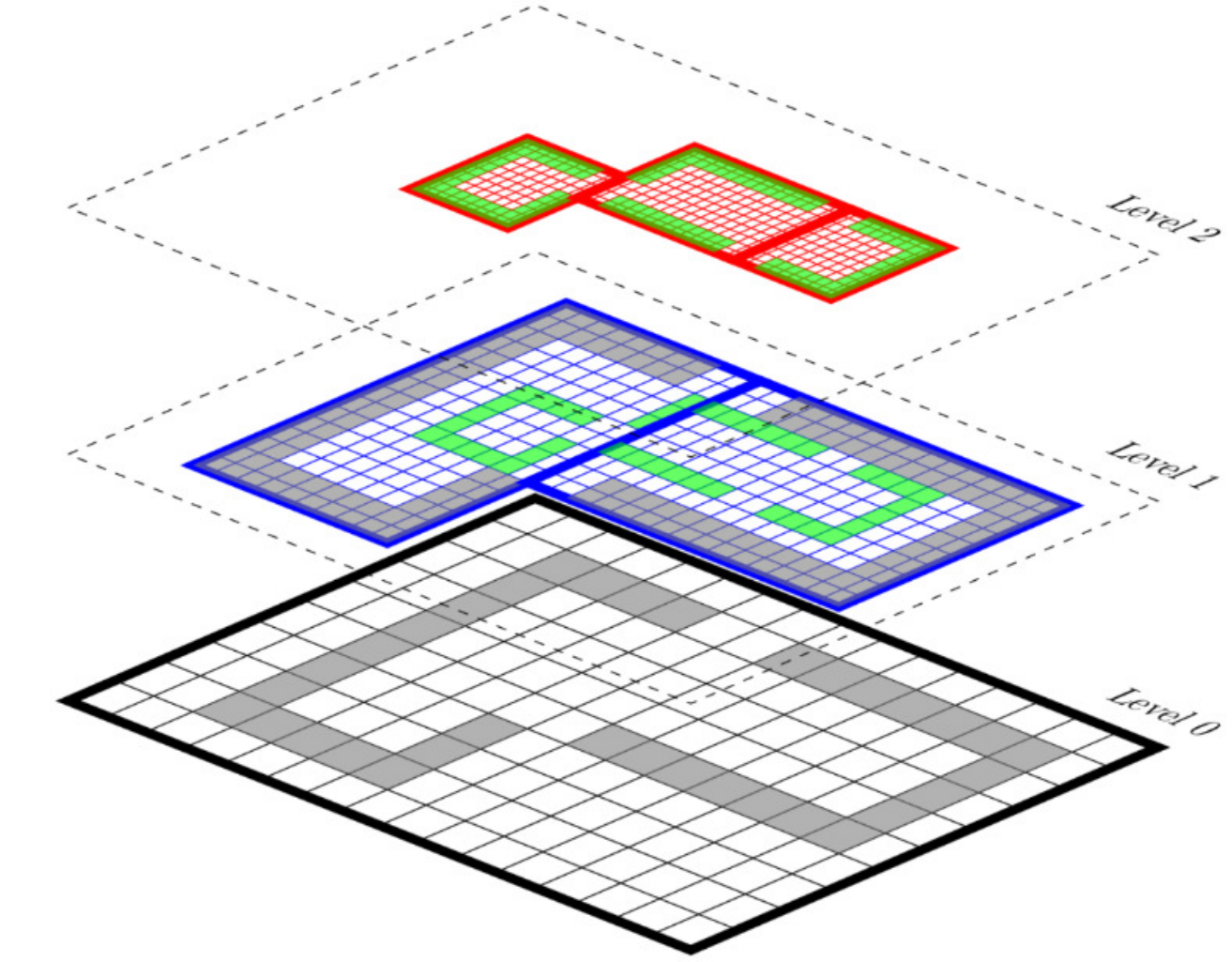}
\caption{Multiple levels in a block-structured AMR grid hierarchy
\label{fig:amr}}
\end{center}
\end{figure}

In block-structured AMR, the solution to a PDE is constructed 
on a hierarchy of levels with different resolution. 
Each level is composed of a union of logically rectangular grids or patches. 
These patches can change dynamically throughout a simulation. In general, 
patches need not be fixed size, and may not have one unique parent grid. 
Figure~\ref{fig:amr} illustrates the use of AMR in a block-structured environment. 


The approach presented by Berger and Oliger~\cite{amr_orig} and
Berger and Colella \cite{Berger} has set the foundation on the patch-based SAMR.
An alternative to the patch-based formulation 
is the octree-based approach which has evolved into the fully-threaded tree (FTT) formalism 
(or cell-based) of Khokhlov \cite{khokhlov1998fully} and the block-based octree of 
MacNiece \textit{et al.} \cite{macneice2000paramesh} \& 
van der Holst \textit{et al.} \cite{van2007hybrid}.

Such AMR methods have gained popularity over the past 30 years and have been 
adopted by various codes in astrophysics.
Some of the well-known examples implementing the patch-based AMR include
AstroBEAR \cite{cunningham2009simulating},
ENZO \cite{bryan2014enzo},
ORION \cite{klein1999star}, 
PLUTO \cite{mignone2011pluto},
CHARM \cite{miniati2011constrained},
CASTRO \cite{almgren2010castro},
MAESTRO \cite{nonaka2010maestro};
the octree-based AMR has been implemented in
FLASH \cite{fryxell2000flash,dubey2009extensible},
NIRVANA \cite{ziegler2008nirvana},
BATS-R-US \cite{powell1999solution,glocer2009multifluid};
the FTT AMR in
RAMSES \cite{teyssier2002cosmological},
ART \cite{kravtsov1997adaptive}.
The AMRVAC code \cite{keppens2003adaptive} features 
both the patch-based and octree-based AMR schemes.

In contrast to these codes 
that incorporate AMR with the purpose of delivering specific applications
in astrophysics, other frameworks have pursued a more general functionality.
Examples include PARAMESH \cite{macneice2000paramesh} 
which supplies solely the octree-based block-structured
mesh capability independent of any governing equations;
AMReX~\cite{amrex} is another a standalone grid software library
that provides the patch-based SAMR support;
Chombo \cite{colella2009chombo,adams2015chombo}
and SAMRAI \cite{hornung2002managing},
on the other hand,
supply both AMR capabilities and a more broader support for solving general systems of equations
of hyperbolic, parabolic, and elliptic partial differential equations (PDEs).
A more compressive survey on the block-structured AMR frameworks
can be found in \cite{dubey2014survey}.

Recently, there have been many noticeable 
efforts aimed at designing high-order accurate solvers
for governing systems of equations (e.g.,
\cite{ray2007using,
mccorquodale2015adaptive,zhang2011order,
buchmuller2014improved,mccorquodale2011high,
dumbser2013ader,
balsara_higher-order_2017,
shu2016high,
reyes_new_2016,
reyes2019variable})
in accordance with a trend of
decreasing memory per compute core in newer
high-performance computing (HPC) architecture
\cite{Attig2011,Dongarra2012future,Subcommittee2014top}.
Such high-order (4th or higher) PDE solvers are then combined with 
the AMR strategies described above.
%
%

Traditionally, a second-order linear interpolation scheme has been commonly adopted
for data prolongation from coarse to finer AMR levels, and 
a mass-conserving averaging scheme for data restriction from finer to coarser levels.
This ``low-order'' AMR interpolation model
has been the default choice in the vast majority of the aforementioned 
AMR paradigms and algorithms in practice.
The accuracy gap between the underlying high-order PDE solvers
and the second-order AMR interpolation could potentially degrade the quality of
solutions from the high-order PDE solvers 
when the solutions are projected to AMR grids that are progressively
undergoing refinements and de-refinements.
In addition, another accuracy loss inevitably happens at fine-coarse boundaries.
It is therefore natural to close the accuracy gap in the direction of providing
high-order models in AMR interpolations, to serves better to maintain the overall
solution accuracy integrated as a whole on AMR grid configurations.
The high-order AMR prolongations of Shen \textit{et al.} \cite{shen2011adaptive} 
and Chen \textit{et al.} \cite{chen_5th-order_2016} are in this vein.
These authors coupled high-order finite difference method (FDM) PDE solvers with 
fourth- or fifth-order accurate prolongations based on the well-known
high-order polynomial interpolation schemes of 
WENO \cite{sebastian2003multidomain} and MP5 \cite{suresh1997accurate}, respectively.
These studies have shown that the AMR simulations with a higher-order coupling
can produce better results in terms of increasing solution accuracy and
lowering numerical diffusion, thereby, resolving fine-scale flow features.

The present work focuses on developing a new high-order polynomial-free 
interpolation scheme
for AMR data prolongation on the block-structured AMR implementation using
the AMReX library. 
Our high-order prolongation scheme stems from the previous studies on 
applying Gaussian Process Modeling \cite{rasmussen2005} in
designing high-order reconstruction/interpolation  
in finite volume method (FVM) \cite{reyes2019variable}
and in finite difference method (FDM)  \cite{reyes2019variable}.

This paper is organized as follows. In Section \ref{sec:AMReX} we overview the relevant
AMR framework, AMReX, as our computational toolkit in which we integrate our new GP-based
prolongation algorithm. In Section \ref{sec:GP} we provide a mathematical overview on
the GP modeling specific to high-order AMR prolongation.
We give step-by-step execution details of our algorithm in Section \ref{sec:method}.
Also, we give a description on extending our work to a GPU-friendly implementation 
by following AMReX programming directives.
Section \ref{sec:results} shows the code performance of the new GP prolongation on
selected multidimensional test problems, and finally, in Section \ref{sec:conclusion}
we summarize the main results of our work.

\section{Overview of AMReX}
\label{sec:AMReX}
Developed and managed by the Center for Computational Science and Engineering at 
Lawrence Berkeley National Laboratory, AMReX is funded through the Exascale Computing Project (ECP) as a software
framework to support the development of block-structured AMR applications focusing on current and next-generation 
architectures~\cite{amrex}. 
AMReX provides support for many operations involving adaptive meshes including multilevel synchronization 
operations, particle and particle/mesh algorithms, solution of parabolic and elliptic systems using geometric and 
algebraic multigrid solvers, and explicit/implicit mesh operations. 
As part of an ECP funded project, AMReX takes the hybrid
MPI/OpenMP CPU parallelization along with GPU implementations (CUDA). 
AMReX is mostly comprised of source files that are written in C++ and Fortran. Fortran is solely 
used for mathematics drivers, while C++ is used for I/O, flow control, memory management and mathematics 
drivers. 

The novelty of the current study is the new GP-based prolongation method
implemented within the AMReX framework. 
The GP implementation furnishes 
an optional high-order prolongation method from coarse to fine AMR levels, 
alternative to the default second-order linear prolongation method in AMReX.
In this way, the GP results in Section \ref{sec:results} naturally inherit all the generic AMReX operations 
such as load balancing, guardcell exchanges, refluxing, AMR data and grid managements,
except for the new GP prolongation method.
We display a suite of test comparisons between the two prolongation methods.

AMR restriction is another important operation on the AMR data management in the opposite direction, 
from fine to coarse levels. We use the default restriction method of averaging
that maintains conservation on AMR grid hierarchies.
This approach populates data on coarse levels 
by averaging down the corresponding fine level data according to
\begin{equation}
\mathbf{U}^C = \frac{1}{R}\sum\limits_i^R{\mathbf{U}^f_i},
\label{eq:avg}
\end{equation}
where $\mathbf{U}^C$ and $\mathbf{U}^f$ are conservative quantities 
on the coarse and fine grids respectively, 
$R=\prod\limits_d{r_d}$ is the normalization factor with  $r_d$ being the 
refinement ratio in each direction $d= x, y, z$. 

Lastly, maintaining conservation across fine-coarse interface levels is done by the operation called the refluxing. 
This process corrects the coarse grid fluxes by averaging down
the fluxes computed on the fine grids abutting the coarse grid.
In practice, the conservation is managed as a posterior correction step 
after all fluid variables ${\bf{U}}^C$ on a coarse cell are updated.
For other AMR operations related to AMReX, interested readers are encouraged to refer to
\cite{zhang2016boxlib,zingale2018meeting,amrex}.

\section{Gaussian Process Modeling for CFD}
\label{sec:GP}
The new prolongation method we are presenting in this paper is based on Gaussian Process (GP) Modeling. 
In order that this paper to be self-contained we give a brief overview on constructing a GP Model 
in this section. More detailed introductions to GP modeling are found in \cite{rasmussen2005,bishop2007pattern}.

Gaussian Processes are a family of stochastic processes in which any finite collection of random variables sampled 
from this process are joint normally distributed. In a more general sense, GPs take samples of functions
from an infinite dimensional function space. In this way, the AMR prolongation
routine described in detail in Section~\ref{sec:method} 
will be drawn from a data-informed distribution space trained on the coarse grid data. 

\subsection{A statistical introduction to Gaussian Processes}
\label{sec:GP_intro}
The construction of the posterior probability distribution over the function space is the heart of 
GP modeling. To construct a GP, one needs to specify a \textit{prior probability distribution} for
the function space. This can be done by specifying two functions, a prior mean function 
and a prior covariance kernel function (see more details below),
by which a GP is fully defined.
Samples, namely function values evaluated at known locations, drawn from the GP prior
are then used to further update this prior probability distribution. 
As a consequence,
a \textit{posterior probability distribution} is generated 
as a combination of the newly updated prior along with these samples, by means of Bayes' Theorem. 

Once constructed, one can draw functions from this data-adjusted GP posterior space 
to generate a model for prolongation in AMR or interpolation/reconstruction.
Specifically, the GP posterior could be used to probabilistically predict the value of a function at points where the 
function has not been previously sampled.
In~\cite{reyes2019variable, reyes_new_2016} Reyes \textit{et al.} have utilized this 
\textit{posterior mean function} as a high-order
predictor to introduce a new class of high-order reconstruction/interpolation algorithms for solving systems of
hyperbolic equations. Within a single algorithmic framework, 
the new GP algorithms have shown a novel algorithmic flexibility in which
a variable order of spatial accuracy is achieved and is given by $2R+1$, corresponding to the size of the one dimensional stencil. 
Here, $R$ is the radius of a GP stencil, called the GP radius, 
given as a positive integer value which represents 
the radial distance between the central cell $\mathbf{x}_i$ and $\mathbf{x}_{i+R}$.

Similarly, from the perspective of designing a probabilistically driven prediction
of function values, the posterior mean function becomes an 
AMR prolongator that delivers a high-order accurate approximation
at the desired location in a computational domain.

As briefly mentioned, GPs can be fully defined by two functions:
a mean function $\bar{f}(\mathbf{x}) = \mathbb{E}[f(\mathbf{x})]$ and a  
covariance function which is a symmetric, positive-definite kernel
$K(\mathbf{x}, \mathbf{y}): \mathbb{R}^N\times\mathbb{R}^N \to \mathbb{R}$. 
Notationally, we write $f\sim \mathcal{GP}(\bar{f}, K)$ to
denote that functions $f$ have been distributed in accordance with
the mean function $\bar{f}(\mathbf{x})$ and the covariance 
$K(\mathbf{x}, \mathbf{y})$ of the GP prior.
Analogous to finite-dimensional 
distributions we write the covariance as 
\begin{equation} 
K(\mathbf{x}, \mathbf{y}) = \mathbb{E}\left[\left(f(\mathbf{x}) - \bar{f}(\mathbf{x})\right)
\left(f(\mathbf{y}) - \bar{f}(\mathbf{y})\right)\right]
\end{equation}
where $\mathbb{E}$ is with respect to the GP distribution. 

One controls the GP by specifying both $\bar{f}(\mathbf{x})$ and $K(\mathbf{x}, \mathbf{y})$, 
typically as some functions parametrized by the so-called hyperparameters. 
These hyperparameters allow us to give the 
``character" of functions generated by the posterior (i.e., length scales, differentiability or regularity)
which will define the underlying pattern of predictions using the posterior GP model.
Suppose we have a given GP and $N$ locations, $\mathbf{x}_n\in \mathbb{R}^d$, 
where $d = 1, 2, 3$ and $ n = 1, \dots, N$. 
For samples $f(\mathbf{x}_n)$ collected at those points, we can calculate the likelihood $\mathcal{L}$, viz., 
the probability of the data $f(\mathbf{x}_n)$ given the GP model. 
Let us denote the data array in a compact form, 
$\mathbf{f} = \left[f(\mathbf{x}_1), \dots, f(\mathbf{x}_N) \right]^T $. 
The likelihood $\mathcal{L}$ of $\mathbf{f}$ is given by
\begin{equation} 
\mathcal{L} \equiv P\big(\mathbf{f} | \mathcal{GP}(\bar{f}, K)\big) = (2\pi)^{-N/2} \det |\mathbf{K}|^{-1/2} 
\exp\left[-\frac{1}{2}\left(\mathbf{f} - \bar{\mathbf{f}}\right)\mathbf{K}
\left(\mathbf{f} - \bar{\mathbf{f}}\right)\right],
\label{eq:likely}
\end{equation} 
where $\mathbf{K}$  is a matrix generated by 
$K_{n,m} = K(\mathbf{x}_n, \mathbf{x}_m)$, $n, m = 1,\dots, N$, 
and the mean $\bar{\mathbf{f}} = [\bar{f}(\mathbf{x}_1), \cdots \bar{f}(\mathbf{x}_N)]^T$. 
Since these samples (or functions) are probabilistically distributed according to the GP prior, i.e., 
$f \sim \mathcal{GP}(\bar{f}, K)$,
we now can make a probabilistic statement about the value of 
any agnostic function $f$ in the GP at a new point $\mathbf{x}_*$,
at which we do not know the exact function value, $f(\bx_*)$.
In other words, the GP model enables us to predict the value of $f(\mathbf{x}_*)$ probabilistically based on the character of likely functions given in the GP model prior. 
For AMR, this is especially important as we need to construct data at a finer resolution where
we do not know the data values at newly generated grid locations refined from a parent coarse level.

An application of Bayes' Theorem, along with the conditioning property,
directly onto the joint Gaussian prior gives the updated (or data-informed) 
posterior distribution of the predicted value $f_*$ conditioned on the observations $\mathbf{f}$, 
\begin{equation} 
P(f_* | \mathbf{f}) = (2\pi U^2)^{-1/2} \exp\left[- \frac{(f_* - \tilde{f}_*)^2}{2U^2}\right],
\label{eq:cond_like}
\end{equation}
where $\tilde{f}_*$ is the posterior mean, given as
\begin{equation}
\tilde{f}_* \equiv \bar{f}(\mathbf{x}_*) + \mathbf{k}_*^T\mathbf{K}^{-1}  (\mathbf{f} - \bar{\mathbf{f}}),
\label{eq:mean}
\end{equation}  
and the \textit{posterior covariance function} as
\begin{equation} 
U^2 \equiv k_{**} - \mathbf{k}_*^T\mathbf{K}^{-1}  \mathbf{k}_*.
\label{eq:cov}
\end{equation} 

The posterior probability given in Eq.~(\ref{eq:cond_like}) is maximized by the choice $f_*=\tilde{f}_*$, leading to Eq.~(\ref{eq:mean}) being taken as the GP prediction for the unknown $f(\bx_*)$. Meanwhile the \textit{posterior} covariance in Eq.~(\ref{eq:cov}) reflects the GP model's confidence in the prediction for the function at $\bfx_*$. we then focus on 
the posterior mean which will become the basis for our interpolation in the GP-based 
AMR prolongation. 

In the next subsections we describe two modeling schemes for the GP-based AMR prolongation. 
In Section~\ref{sec:gp_ptwise_prol} we describe the first method that prolongates pointwise state data
from coarse to fine levels. 
For AMR simulations in which the state data is represented as volume-averaged, conserving
such quantities become crucial to satisfy the underlying conservation laws. 
To meet this end, 
we introduce the second method in Section~\ref{sec:gp_vol_prol}, 
which preserves volume-averaged quantities in prolongation.
We will refer to our GP-based AMR prolongation as GP-AMR for the rest of this paper.


\subsection{GP for pointwise AMR prolongation}
\label{sec:gp_ptwise_prol}
In this section we introduce the first GP-AMR prolongation method that is suitable for
AMR applications where the state data is comprised of pointwise values. 
In this case the GP-AMR model samples are given as pointwise evaluations of the underlying function.
Let $\Delta x_d$ denote the distance between points in a
\textbf{coarse} level in each $d=x,y,z$ direction.  
Using the posterior mean function in Eq.~\eqref{eq:mean},
we first devise a pointwise prolongation scheme for 
AMR, i.e., AMR prolongation of pointwise data from coarse to fine levels.  
The choice of $\mathbf{x}_*$ will depend on the \textit{refinement ratio} 
$\mathbf{r} = [r_x, r_y, r_z]$ and there will 
be $\prod_d r_d $ new points generated for the new level in general. 
For example, if we wished to refine a single coarse grid by two in 
all three directions in 3D, we would generate eight new grid points as well as 
the eight new  associated data values at those grid points
in a newly refined level. 
To illustrate the process, we consider a simple example of a two-level refinement in 1D. 
In this refinement two refined data values are to be newly generated for each and every coarse value. 
Assume here that we utilize a stencil with the GP radius of one (i.e., $R=1$) in which case the local 3-point 
GP stencil $\mathbf{f}_i$ centered at each $i$-th cell for interpolation is laid out as
\[\mathbf{f}_i = [q_{i-1}, q_i, q_{i+1}]^T. \]
From this given stencil data at the coarse level, 
we wish to generate two finer data values $q_{s\pm1/2}$
for each $s=i-1, i, i+1$. 
To do this, we use the posterior mean function in Eq.~\eqref{eq:mean} 
on three 3-point GP stencils, $\mathbf{f}_s$, $s=i-1, i, i+1$, 
to populate a total of six new data,
\begin{equation}
q_{s\pm\frac{1}{2}} = \mathbf{k}_{s\pm\frac{1}{2}}^T\mathbf{K}^{-1} \mathbf{f}_s, \;\; s=i-1, i, i+1,
\label{eq:pinterp} 
\end{equation}
where we used a zero mean prior, $\bar{\mathbf{f}}=0$.
In 2D or 3D, data values on a standard $(2N+1)$-point stencil 
are to be reshaped into a 1D local array $\mathbf{f}_s$
in an orderly fashion, where each $\mathbf{f}_s$ includes corresponding multidimensional 
data reordered in 1D between $s-N$ and $s+N$.
This strategy will be fully described in Section~\ref{sec:method}. 

A common practice with GP modeling is to assume 
a zero prior mean as we did with Eq.~\eqref{eq:pinterp}. 
In our implementations we use this assumption.
Something to note is that the GP weights, $\mathbf{k}_{*}^T\mathbf{K}^{-1}$,
are independent of the 
samples $\bf{f}$, and are constructed based on 
the choice of kernel function and the location of the samples, $\bf{x}_n$, 
and prediction point, $\bf{x}_*$, alone. 
This is particularly useful in block structured AMR applications, as we can compute the weights for 
each level a priori, based on the min and max levels prescribed for each run. 
Otherwise, we can generate the model weights the first time a level is used and save them for later uses.

Since the matrix $\mathbf{K}$ is symmetric and positive-definite, 
we can use the Cholesky decomposition to compute the GP weights. 
In practice, 
we compute and save $\bf{w}^T_* = \bf{k}^T_*\bf{K}^{-1}$ using Cholesky followed
by back-substitution only once per simulation,
either at an initial grid configuration step or at the first time an AMR level is newly used.
In this way,
the computational cost of the prolongation is reduced to
a dot product between $\bf{w}$ and $\bf{f}$.
As a consequence we arrive at a compact form,
\begin{equation}
q_{s\pm\frac{1}{2}} = \mathbf{w}^T_{s\pm\frac{1}{2}} \mathbf{f}_s, \;\; s=i-1, i, i+1.
\label{eq:pinterp2} 
\end{equation}

There are many choices of covariance kernel functions available for GP modelling \cite{rasmussen2005,bishop2007pattern}. One of the most
widely used kernels in Gaussian Process modeling is 
the squared-exponential (SE) covariance kernel function,
\begin{equation}
K(\mathbf{x}, \mathbf{y}) \equiv \Sigma^2\exp\left[-\frac{(\mathbf{x} - \mathbf{y})^2}{2\ell^2}\right].
\label{eq:sqrexpcov}
\end{equation}
The SE kernel is infinitely differentiable and as a consequence will sample functions that are then equally smooth. The kernel contains two model hyperparameters $\Sigma$ and $\ell$. $\Sigma$ acts as an overall constant factor that has no impact on the posterior mean (this can be seen as a cancellation between the $\bfk_*^T$ and $\bfK^{-1}$ terms in Eq.~(\ref{eq:mean})), which is the basis of our GP-AMR prolongation and we take $\Sigma=1$. $\ell$ controls the length scale on which \textit{likely} functions will vary according to the GP model. 

All that remains to complete the GP model is to specify the prior mean function. The prior mean function is often depicted as a constant mean function for simplicity, i.e., 
$\bar{f}(\mathbf{x}) = f_0 $. $f_0$ controls the behavior of the GP prediction at spatial locations that, according to the kernel function, are not highly correlated with any of the observed values $\bff$. In the context of the SE kernel for prolongation this happens when the choice of $\ell$ is much smaller than the grid spacing between coarse cells. For that reason it is advisable to choose $\ell$ so that it is on the order of the size of the prolongation stencil.
%
%
The model using the SE kernel in Eq.~\eqref{eq:sqrexpcov} and Eq.~(\ref{eq:pinterp}) with the prescribed
hyperparameter choices is our first formula for the pointwise AMR prolongation.


\subsection{A GP Prolongation for Cell-Averaged Quantities}
\label{sec:gp_vol_prol}
For the majority of AMReX and fluid dynamics application codes, 
the state data is cell-averaged (or volume-averaged),
as per the formulation of FVMs.
The above GP-AMR prolongation for pointwise data has to be modified in order to preserve the integral relations between fine and coarse data that are implicit in the integral formulation of the governing equations used in FVM.
The key observation from \cite{reyes_new_2016} is that the averaging over cells constitutes a "linear" operation on the a function $f(\bfx)$. As is done for finite dimensional Gaussian processes, linear operations on Gaussian random variables yields a new Gaussian random variable with linearly transformed mean and covariance functions.
In order to 
calculate the covariance between cell averaged quantities 
we need an integrated covariance kernel as described in
~\cite{reyes_new_2016}. 
That is, 
\begin{equation}
\begin{aligned} 
C_{kh} & =  \mathbb{E}[(G_k - \bar{G}_k)(G_h - \bar{G}_h)]  \\
&  =  \int \mathbb{E}[(f(\mathbf{x}) - \bar{f}(\mathbf{x}))(f(\mathbf{y}) - \bar{f}(\mathbf{y}))]dg_k(\mathbf{x})
dg_h(\mathbf{y}) \\ & = 
\iint K(\mathbf{x}, \mathbf{y}) dg_k(\mathbf{x}) dg_h(\mathbf{y}),
\end{aligned}
\label{eq:int}
\end{equation} 
where 
\begin{equation} 
dg_s(\mathbf{x}) = \begin{cases} \displaystyle d\mathbf{x}\prod_{d=x,y,z}^{D} \frac{1}{\Delta x_d} & 
\textrm{if   } \mathbf{x}\in I_s, \\ 
0 & \textrm{else,}\end{cases}
\end{equation}
is the $1D$ cell-average measure. $G_i = \left<f(\mathbf{x}_i)\right> = \frac{1}{\mathcal{V}}\int_{I_i} f(\mathbf{x}_i) d\mathcal{V}$ are the cell-averaged data over cell $I_i\subset \mathbb{R}^D$ with volume $\mathcal{V}$.

With the use of the squared-exponential kernel, Eq.~(\ref{eq:int}) becomes 
\begin{equation}
\begin{aligned}
C_{kh} = \prod_{d=x,y,z}^{D} \sqrt{\pi}\left(\frac{\ell}{\Delta x_d}\right)^2\Bigg\{\left(
\frac{\Delta_{kh}+1}{\sqrt{2}\ell/\Delta x_d}\textrm{erf}\left[\frac{\Delta_{kh}+1}{\sqrt{2}\ell/\Delta x_d}\right]
+ \frac{\Delta_{kh}-1}{\sqrt{2}\ell/\Delta x_d}\textrm{erf}\left[\frac{\Delta_{kh}-1}{\sqrt{2}\ell/\Delta x_d}\right]
\right) \\ 
+ \frac{1}{\sqrt{\pi}}\left(\exp{\left[-\left(\frac{\Delta_{kh}+1}{\sqrt{2}\ell/\Delta x_d}\right)^2\right]}
+\exp{\left[-\left(\frac{\Delta_{kh}-1}{\sqrt{2}\ell/\Delta x_d}\right)^2\right]}\right) 
\\ - 2\left(\frac{\Delta_{kh}}{\sqrt{2}\ell/\Delta x_d}\textrm{erf}\left[\frac{\Delta_{kh}}{\sqrt{2}\ell/\Delta x_d}
\right] + \frac{1}{\sqrt{\pi}}\exp\left[-\left(\frac{\Delta_{kh}}{\sqrt{2}\ell/\Delta x_d}\right)^2\right]\right)
\Bigg\},
\end{aligned}
\label{eq:intsq}
\end{equation}
where we used $\Delta_{kh} = (x_{d,h} - x_{d,k})/\Delta x_d$.  

Following similar arguments as for the covariance kernel function, the prediction vector $\bfk_*$ must also be linearly transformed to reflect the relationship between the input data averaged over coarse cells and the output prolonged data averaged over the fine cells. This leads to
\[T_{k*} \equiv T(\mathbf{x}, \mathbf{x}_*) = \int_{I_k} \int_{I_*} 
K(\mathbf{x}, \mathbf{x}_*) dg_k(\mathbf{x}) dg_*(\mathbf{x_*}), \] 
where 
\[I_* = \bigtimes\limits_{d=x,y,z}^{D} \left[x_{*,d} - \frac{\Delta x_d}{2r_d}, \;\; x_{*,d} + \frac{\Delta x_d}{2r_d}\right],\]
in which $\bigtimes$ denotes the Cartesian production on sets. 
Using the SE kernel, we have a closed form for $T_{k*}$,
\begin{equation} 
T_{k*} = \pi^{D/2}\prod_{d = x,y,z}^D r_d \left(\frac{\ell}{\Delta x_d}\right)^2 \sum_{\alpha = 1}^4 (-1)^\alpha 
\left[\phi_{\alpha, d} \textrm{erf}(\phi_{\alpha, d}) + \frac{1}{\sqrt{\pi}}\exp(-\phi_{\alpha,d}^2)\right],
\label{eq:trans}
\end{equation}  
where for each $\alpha=1,\dots,4$,
\[\phi_{\alpha,d}= \frac{1}{\sqrt{2}\ell/\Delta x_d}\left(\Delta_{k,*} + \frac{r_d -1}{2r_d}, \hspace{1mm}  
	      \Delta_{k,*} + \frac{r_d+1}{2r_d}, \hspace{1mm} \Delta_{k,*} - \frac{r_d -1}{2r_d}, \hspace{1mm}
	      \Delta_{k,*} - \frac{r_d+1}{2r_d}\right).\]
Therefore, with the combination of the cell-averaged kernel in Eq.~(\ref{eq:intsq}) and 
the weight vector in Eq.~(\ref{eq:trans}), we obtain
our second GP-AMR formula given in the integral analog of Eq.~(\ref{eq:mean}) 
for cell-averaged data prolongation from coarse to fine levels,
\begin{equation}
\left<f(\mathbf{x}_*)\right> = \mathbf{T}_{*}^T\mathbf{C}^{-1} \mathbf{G},
\label{eq:vinterp} 
\end{equation}
where we used the zero mean as before. 
The vector $\mathbf{G}$ of cell-averaged samples within the GP radius $R$
is given as $\mathbf{G} = [G_{i-R}, \dots, G_{i+R}]^T$.
Analogous to the pointwise method, we cast $\mathbf{T}_{*}^T\mathbf{C}^{-1}$
into a new GP weight vector $\mathbf{z}_*$ to rewrite Eq.~\eqref{eq:vinterp} as
\begin{equation}
\left<f(\mathbf{x}_*)\right> = \mathbf{z}_{*}^T\mathbf{G}.
\label{eq:vinterp2} 
\end{equation}

Many methods perform interpolation in a dimension-by-dimension manner. 
In contrast, the above two GP-AMR methods are inherently multidimensional. 
Moreover, the use of the SE kernel as a base in each $d$-direction 
facilitates the analytic multidimensional form obtained in Eq.~\eqref{eq:intsq}.
Our GP-AMR methods, therefore, provide a unique framework where all interpolation procedures
in AMR grid hierarchies naturally support multidimensionality, as the evaluation of the covariance matrices only depends on the distance between data points.
Furthermore, it is worth to point out that the two prolongation schemes in 
Eqs.~\eqref{eq:pinterp2} and \eqref{eq:vinterp2} 
are merely a straightforward calculation of
dot products between the GP weight vectors and the grid data.

This is the novelty of the use of GP modeling in AMR prolongation, which reveals two new
compact prolongation methods that are computed in the same way for any stencil configuration in any number of spatial dimensions without any added complexity. This is in stark contrast to polynomial based methods which require the use of explicit basis functions and have strict requirements on stencil sizes and configurations in order to form a well posed interpolation. All of this together results in polynomial based methods giving increased difficulty especially as the number of dimensions or the order of accuracy is increased, highlighting the simplicity afforded in the present GP-AMR methods.


\subsection{Nonlinear Multi-substencil Method of GP-WENO for Non-Smooth Data}
\label{sec:GP-WENO}
Both of the above GP modeling techniques can suffer from non-physical
oscillations 
near discontinuities. 
The SE and integrated SE kernels work very well for continuous data, but we need to implement some type of 
a ``limiting" process in order to suppress `nonphysical' oscillations in flow regions with sharp gradients. The linear interpolations used by default in AMReX makes use of the monotonized central (MC) slope limiter in order to produce slopes that do not introduce any new extrema into the solution.
In this study we utilize the GP-based smoothness indicator approach studied in
GP-WENO \cite{reyes_new_2016,reyes2019variable}.
Following the Weighted Essentially Non-Oscillaltory (WENO) \cite{jiang1996efficient} approach, GP-WENO adaptively chooses a non-oscillatory stencil by nonlinearly weighing GP predictions trained on a set of substencils according to a GP-based local indicator of smoothness $\beta_m$. 
The smoothness is determined using the GP likelihood to measure the compatibility of the substencil data $\bf_m$ with the smooth SE kernel. 
Effectively $\beta_m$ is an indication of how well the data in $\bf_m$ matches with the GP model assumptions that are encoded in $\mathbf{K}_{m,\sigma}$ (or $\mathbf{C}_{m,\sigma}$).
There two differences between
$\mathbf{K}$ and  $\mathbf{K}_{m,\sigma}$ 
in that
(i) $\mathbf{K} \in \mathbb{R}^{M \times M} $ and  
$\mathbf{K}_{m,\sigma} \in \mathbb{R}^{(2D+1)\times(2D+1)}$, where
$M=2D^2+2D+1$ for each spatial dimension $D=1,2,3$, and 
(ii) the scale-length hyperparameter for $\mathbf{K}_{m,\sigma}$,  $\sigma$ is 
a much smaller length scale 
in accordance with the narrow shock-width spread over a couple of grid spacing.
The same differences hold between $\mathbf{C}$ and  $\mathbf{C}_{m,\sigma}$ as well.
%

The first step in this multi-substencil method of GP-WENO
is to build $2D+1$ substencil data on each substencil $S_m$, $m=1, \dots, 2D+1$. 
The data are combined using linear weights $\gamma_m$ derived from an over-determined linear 
system relating the weights generated by building a GP model on all substencils $S_m$ and the weights
generated from a GP model on a total stencil $S$. The last step is to take the linear weights $\gamma_m$
to define nonlinear weights $\omega_m$ using the GP-based smoothness indicators $\beta_m$~\cite{reyes_new_2016, reyes2019variable}.

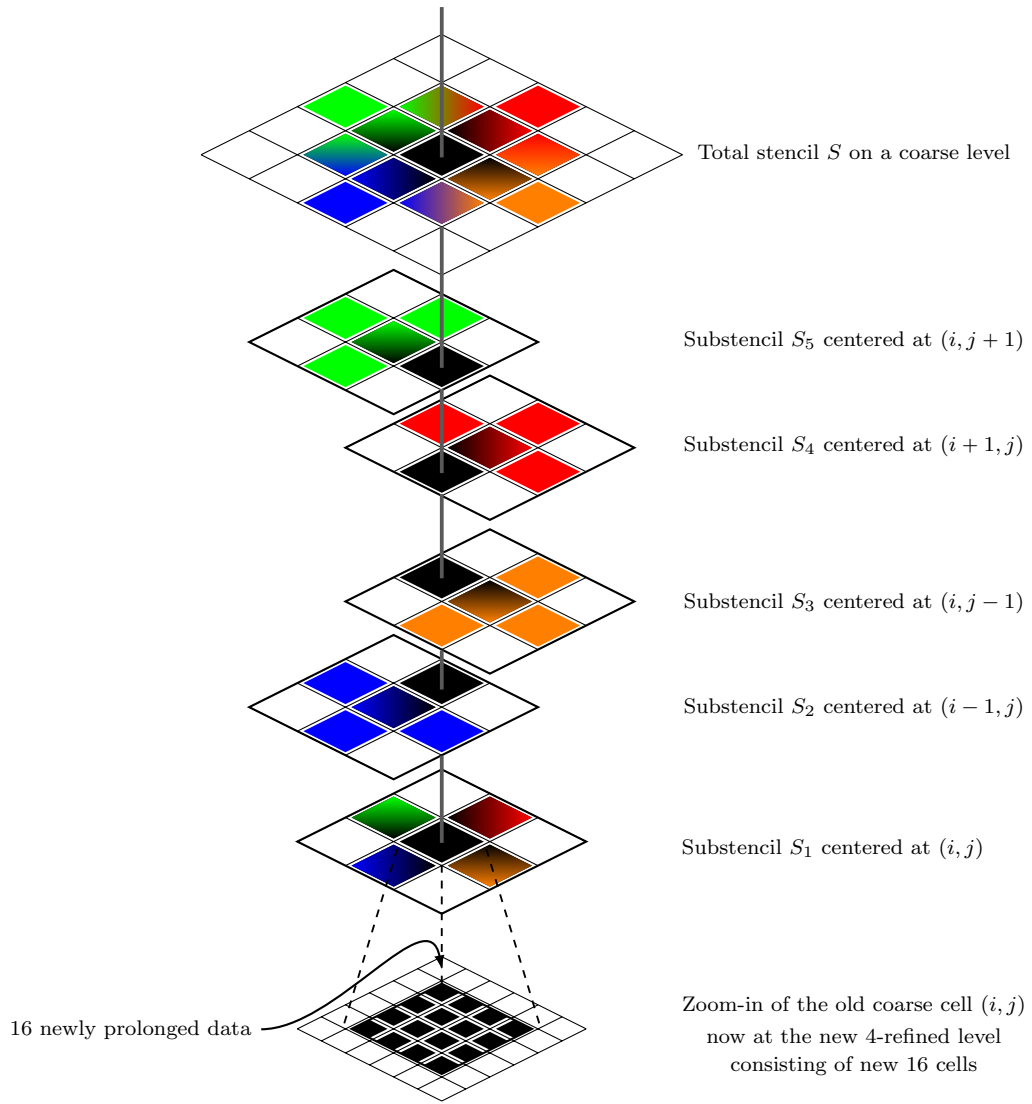
\begin{figure}
\begin{center}
\begin{tikzpicture}[scale=.8,every node/.style={minimum size=1cm}]
\begin{scope}[ yshift=-125,every node/.append style={ yslant=0.5,xslant=-1},yslant=0.5,xslant=-1 ] 
\draw[step=4mm, black] (1.2,1.2) grid (3.6,3.6); 
\draw[black,thick] (1.65,1.65) rectangle (3.15,3.15);
\fill[black] (2.05,2.05) rectangle (2.35,2.35); 
\fill[black] (1.65,2.05) rectangle (1.95,2.35); 
\fill[black] (2.45,2.05) rectangle (2.75,2.35); 
\fill[black] (2.05,2.45) rectangle (2.35,2.75); 
\fill[black] (2.05,1.95) rectangle (2.35,1.65); 
\fill[black] (1.65,2.45) rectangle (1.95,2.75); 
\fill[black] (2.45,2.45) rectangle (2.75,2.75); 
\fill[black] (2.75,1.95) rectangle (2.45,1.65); 
\fill[black] (1.65,1.95) rectangle (1.95,1.65); 
\fill[black] (3.15,1.95) rectangle (2.85,1.65);
\fill[black] (3.15,2.35) rectangle (2.85,2.05);
\fill[black] (3.15,2.75) rectangle (2.85,2.45);
\fill[black] (3.15,3.15) rectangle (2.85,2.85);
\fill[black] (2.75,3.15) rectangle (2.45,2.85);
\fill[black] (2.35,3.15) rectangle (2.05,2.85);
\fill[black] (1.95,3.15) rectangle (1.65,2.85);
\end{scope} %
\begin{scope}[ yshift=150,every node/.append style={ yslant=0.5,xslant=-1},yslant=0.5,xslant=-1 ] 
\fill[white,fill opacity=0.25] (1.6,0.8) rectangle (4,3.2); 
\draw[step=8mm, black] (1.6,0.8) grid (4,3.2); 
\draw[black,thick] (1.6,0.8) rectangle (4,3.2);
\fill[red] (2.45,2.45) rectangle (3.15,3.15); 
\fill[right color =red, left color = black] (3.15,2.35) rectangle (2.45,1.65); 
\fill[black] (1.65,2.35) rectangle (2.35,1.65); 
\fill[red] (3.15,1.55) rectangle (2.45,0.85); 
\fill[red] (3.25, 2.35) rectangle (3.95, 1.65); 
\end{scope}

\begin{scope}[ yshift=200,every node/.append style={ yslant=0.5,xslant=-1},yslant=0.5,xslant=-1 ] 
\fill[white,fill opacity=0.25] (0.8,1.6) rectangle (3.2,4); 
\draw[step=8mm, black] (0.8,1.6) grid (3.2,4); 
\draw[black,thick] (0.8,1.6) rectangle (3.2,4);
\fill[green] (2.45,2.45) rectangle (3.15,3.15); 
\fill[green] (0.85,2.45) rectangle (1.55,3.15); 
\fill[green] (1.65,3.25) rectangle (2.35,3.95); 
\fill[bottom color =black, top color =green] (1.65,2.45) rectangle (2.35,3.15); 
\fill[black] (1.65,2.35) rectangle (2.35,1.65); 
\end{scope}

\begin{scope}[ yshift=-25,every node/.append style={ yslant=0.5,xslant=-1},yslant=0.5,xslant=-1 ] 
\fill[white,fill opacity=0.25] (0.8,0.8) rectangle (3.2,3.2); 
\draw[step=8mm, black] (0.8,0.8) grid (3.2,3.2); 
\draw[black,thick] (0.8,0.8) rectangle (3.2,3.2);
\fill[top color =green, bottom color = black] (1.65,2.45) rectangle (2.35,3.15); 
\fill[left color = black, right color = red] (3.15,2.35) rectangle (2.45,1.65); 
\fill[black] (1.65,2.35) rectangle (2.35,1.65); 
\fill[left color =blue, right color = black] (0.85,2.35) rectangle (1.55,1.65); 
\fill[top color=black, bottom color =orange] (1.65,1.55) rectangle (2.35,0.85); 
\end{scope}

\begin{scope}[ yshift=50,every node/.append style={ yslant=0.5,xslant=-1},yslant=0.5,xslant=-1 ] 
\fill[white,fill opacity=0.25] (0,0.8) rectangle (2.4,3.2); 
\draw[step=8mm, black] (0,0.8) grid (2.4,3.2); 
\draw[black,thick] (0,0.8) rectangle (2.4,3.2);
\fill[blue] (0.85,2.45) rectangle (1.55,3.15); 
\fill[black] (1.65,2.35) rectangle (2.35,1.65); 
\fill[left color =blue, right color = black] (0.85,2.35) rectangle (1.55,1.65); 
\fill[blue] (0.85,1.55) rectangle (1.55,0.85);
\fill[blue] (0.05, 2.35) rectangle (0.75, 1.65);
\end{scope}

\begin{scope}[ yshift=100,every node/.append style={ yslant=0.5,xslant=-1},yslant=0.5,xslant=-1 ] 
\fill[white,fill opacity=0.25] (0.8,0) rectangle (3.2,2.4); 
\draw[step=8mm, black] (0.8,0) grid (3.2,2.4); 
\draw[black,thick] (0.8,0) rectangle (3.2,2.4);

\fill[black] (1.65,2.35) rectangle (2.35,1.65); 
\fill[top color = black, bottom color = orange] (1.65,1.55) rectangle (2.35,0.85); 
\fill[orange] (3.15,1.55) rectangle (2.45,0.85); 
\fill[orange] (0.85,1.55) rectangle (1.55,0.85);
\fill[orange] (1.65,0.05) rectangle (2.35,0.75);
\end{scope}

\begin{scope}[ yshift=300,every node/.append style={ yslant=0.5,xslant=-1},yslant=0.5,xslant=-1 ] 
\fill[white,fill opacity=0.25] (0,0) rectangle (4,4); 
\draw[step=8mm, black] (0.0,0) grid (4,4); 
\fill[green] (1.65,3.25) rectangle (2.35,3.95); 
\fill[left color = green, right color =red] (2.45,2.45) rectangle (3.15,3.15); 
\fill[red] (3.25, 2.35) rectangle (3.95, 1.65); 
\fill[top color = red, bottom color =orange] (3.15,1.55) rectangle (2.45,0.85); 
\fill[orange] (1.65,0.05) rectangle (2.35,0.75);
\fill[top color = green, bottom color =blue] (0.85,2.45) rectangle (1.55,3.15); 
\fill[left color =blue, right color = orange] (0.85,1.55) rectangle (1.55,0.85);
\fill[blue] (0.05, 2.35) rectangle (0.75, 1.65);
\fill[top color =green, bottom color = black] (1.65,2.45) rectangle (2.35,3.15); 
\fill[left color = black, right color = red] (3.15,2.35) rectangle (2.45,1.65); 
\fill[black] (1.65,2.35) rectangle (2.35,1.65); 
\fill[left color =blue, right color = black] (0.85,2.35) rectangle (1.55,1.65); 
\fill[top color=black, bottom color =orange] (1.65,1.55) rectangle (2.35,0.85); 
\end{scope}

\draw[line width=0.5mm, gray!70!black](0, 15) -- (0,12.5);
\draw[line width=0.5mm, gray!70!black](0, 11.35) -- (0,9);
\draw[line width=0.5mm, gray!70!black](0, 8.65) -- (0,7.25);
\draw[line width=0.5mm, gray!70!black](0, 6.9) -- (0,5.5);
\draw[line width=0.5mm, gray!70!black](0, 4.3) -- (0,3.65);
\draw[line width=0.5mm, gray!70!black](0, 2.55) -- (0,1.1);

\draw[line width=0.25mm, dashed](-.74,1) -- (-1.65, -2);
\draw[line width=0.25mm, dashed](0.74,1) -- ( 1.65, -2);
\draw[line width=0.25mm, dashed](0, 1) -- (0, -1.5); 

\draw[-latex,thick,black](-3,-2)node[left]  
{\footnotesize{16 newly prolonged data}} to[out=0,in=100] (-0,-1); 
\draw[thick,black](6.5, 1) node {\footnotesize{Substencil $S_1$ centered at $(i,j)$}};
\draw[thick,black](6.85, 3.35) node {\footnotesize{Substencil $S_2$ centered at $(i-1, j)$}}; 
\draw[thick,black](6.85, 5.1) node {\footnotesize{Substencil $S_3$ centered at $(i, j-1)$}}; 
\draw[thick,black](6.85, 7.7) node {\footnotesize{Substencil $S_4$ centered at $(i+1, j)$}}; 
\draw[thick,black](6.85, 9.45) node {\footnotesize{Substencil $S_5$ centered at $(i, j+1)$}};
\draw[thick,black](6.85, 12.6) node  {\footnotesize{Total stencil $S$} on a coarse level};

\draw[thick,black](6.85,-1.6) node {\footnotesize{Zoom-in of the old coarse cell $(i,j)$}};
\draw[thick,black](6.85,-2.1) node {\footnotesize{now at the new 4-refined level}}; %
\draw[thick,black](6.85,-2.6) node {\footnotesize{consisting of new 16 cells}}; %
\end{tikzpicture}
\end{center}
\caption{GP prolongation using five GP substencils that are combined to 
produce 16 new prolonged data on a 2D finer grid. The 4-refinement ratio in both
$x$ and $y$ directions is considered here to prolong the single data 
from the old coarse cell $(i,j)$ to 16 newly refined locations.
\label{fig:2dS}}
\end{figure}

We now describe in detail the GP-AMR method for the two-dimensional case. 
Extensions to other dimensions are readily obtained due to the choice of the isotropic SE kernel, with the only difference being in the number of stencil points used. 
We begin with a total stencil $S$, taken as all cells whose index centers are within a radius of $2$ of the central cell $I_{i,j}$.
The total stencil $S$ is then subdivided into $2D+1$ candidate stencils, $S_m, m = 1, \dots 2D+1$,  such that 
$\displaystyle\bigcap_{m=1}^{2D+1} S_m = \{\mathbf{x}_{i,j} \}$ 
and $\displaystyle\bigcup_{m=1}^{2D+1} S_m = S$. A schematic of these stencil configurations is given in Fig.~\ref{fig:2dS}.
That is, the prolongation will have the form: 
\begin{equation}
f_* = \sum_{m=1}^{2D+1} \omega_m \mathbf{w}_m^T\mathbf{f}_m  
\label{eq:MSGP} 
\end{equation} 
where $\mathbf{w}_m^T = \mathbf{T}_{*,m}\mathbf{C}^{-1}_{m}$ for the cell-averaged prolongation, or 
$\mathbf{w}_m^T = \mathbf{k}_{*,m}^T\mathbf{K}^{-1}_{m}$ for the pointwise prolongation. 
The coefficients $\omega_m$ are defined as in the WENO-JS method~\cite{WENO},  
\beq
\omega_m = \frac{\tilde{\omega}_m}{\sum_s \tilde{\omega}_s} \quad \textrm{where} \quad 
\tilde{\omega}_m = \frac{\gamma_m}{(\epsilon + \beta_m)^p}.
\eeq
For our algorithm we choose  $\epsilon = 10^{-36}$ and $p = 2$. The terms $\beta_m$ are 
taken as the data dependent term in $\log(\mathcal{L})$ from Eq.~(\ref{eq:likely}) for the stencil data
$\mathbf{f}_m$, that is, 
\beq
\beta_m = \mathbf{f}_m^T \mathbf{K}_{m,\sigma}^{-1} \mathbf{f}_m
\eeq
for the pointwise prolongation, and 
\beq
\beta_m = \mathbf{f}_m^T \mathbf{C}_{m,\sigma}^{-1} \mathbf{f}_m
\eeq
for the cell-averaged prolongation. 
Notice that, due to the properties of the 
kernel matrices~\cite{reyes_new_2016}, we can cast 
\begin{equation} 
\beta_m = \sum_{i=1}^{2D+1} \frac{1}{\lambda_i}\left(\mathbf{v}_i^T\mathbf{f}_m\right)^2,
\label{eq:eig} 
\end{equation} 
where $\mathbf{v}_i$ and $\lambda_i$ are the eigenvectors and eigenvalues of the covariance kernel
matrix, ${\bf{K}}_{m,\sigma}$ or ${\bf{C}}_{m,\sigma}$.
As described in~\cite{reyes_new_2016,reyes2019variable}, the GP-based smoothness indicators $\beta_m$ defined in this way
is derived by taking the negative log of the GP likelihood of Eq.~\eqref{eq:likely}.
This gives rise to the statistical interpretation of $\beta_m$ which relates that
if there is a shock or discontinuity in one of the substencils, say $S_k$, such a short length-scale (or rapid) 
change on $S_k$ makes $\mathbf{f}_k$ unlikely. Here this likeliness is relative to a GP model that assumes that the underlying function is smooth on the length scale set by $\sigma$. In other words, the GP model whose smoothness is
represented by the smoothness property of its covariance kernel, $\mathbf{K}_{m,\sigma}$ or 
$\mathbf{C}_{m,\sigma}$,
gives a low probability to  $\mathbf{f}_k$, in which case $\beta_k$ -- given as the negative log 
likelihood of $\mathbf{f}_k$ -- becomes relatively larger than the other $\beta_m$, $m \ne k$.

In our method we use GP modeling for both a regression (prolongation) and 
a classification. The regression aspect enables us to prolongate GP samples 
(i.e., function values, or fluid values) over the longer length-scale variability specified by $\ell$.
On the other hand, the classification aspect allows us to detect and handle discontinuities. This
is achieved by employing a much shorter length-scale variability tuned by $\sigma$,
which is integrated into the eigensystem in Eq.~(\ref{eq:eig}) 
generated with $\mathbf{K}_{m,\sigma}$ or $\mathbf{C}_{m,\sigma}$. 
Smaller than $\ell$,
the parameter $\sigma$ is chosen to reflect the short width of shocks and discontinuities in 
numerical simulations, which is typically over a couple of grid spacings.
%
In this manner, we use 
two length scale parameters, $\ell$ for the interpolation model, and $\sigma$ for shock-capturing.  

Another key factor are the linear weights $\gamma_m$, $m=1, \dots, 2D+1$.
Let $\boldsymbol{\gamma}=[\gamma_1, \dots, \gamma_{2D+1}]^T$ be a vector
containing $2D+1$ linear weights, each corresponding to one of the substencils. 
These weights are retrieved by solving an over-determined linear system
\begin{equation}
\mathbf{M}\boldsymbol{\gamma} = \mathbf{w}_*,
\label{eq:gamsys} 
\end{equation}
where the $n$-th column of $\mathbf{M}$ is given by $\mathbf{w}_n$, and $\bf{w}_*$ is the model weights for 
the interpolation point $\bf{x}_*$ relative to the total stencil $S$. As mentioned previously, 
these weights are generated using the length scale parameter $\ell$. We should note that $\bf{M}$ is a potentially sparse 
matrix, and is constructed using the substencil model weights.  

For our GP modeling procedure in multiple spatial dimensions, 
we cast the multidimensional stencil $S$ as a flattened array. 
To illustrate this concept we explore a 2$D$ example where
the coarse level cells are refined by the 4-refinement ratio
in both $x$ and $y$ directions, i.e., $r_x=r_y=4$.
Suppose $D = 2$, 
in which case and the total stencil $S$ is in the $5\times5$ patch of cells centered at $(i,j)$ and
 contains 13 data points. The total stencil
is subdivided into five 5-point substencils
$S_m$, $m=1, \dots, 5$.
We take the natural cross-shape substencil for each $S_m$ 
on each of which GP will
approximate function values (i.e., state values of density, pressure, etc.) 
at 16 new refined locations, i.e.,
$(i\pm 1/4, j \pm 1/4)$,
$(i\pm 1/4, j \pm 3/4)$,
$(i\pm 3/4, j \pm 1/4)$, and
$(i\pm 3/4, j \pm 3/4)$.

For instance, let's choose $\mathbf{x}_{i+1/4, j+1/4}$ as the location
we wish GP to compute function values for prolongation.
Explicitly, five 5-point substencils are chosen as,
\bea
S_1 &=& \left[\mathbf{x}_{i,j-1}, \mathbf{x}_{i-1,j}, \mathbf{x}_{i,j}, \mathbf{x}_{i+1,j}, \mathbf{x}_{i,j+1}\right], \nonumber \\
%
%
S_2 &=& \left[\mathbf{x}_{i,j-2}, \mathbf{x}_{i-1,j-1}, \mathbf{x}_{i,j-1}, \mathbf{x}_{i+1,j-1}, \mathbf{x}_{i,j}\right],\nonumber\\
S_3 &=& \left[\mathbf{x}_{i+1,j-1}, \mathbf{x}_{i,j}, \mathbf{x}_{i+1,j}, \mathbf{x}_{i+2,j}, \mathbf{x}_{i+1,j+1}\right], \\
S_4 &=& \left[\mathbf{x}_{i,j}, \mathbf{x}_{i-1,j+1}, \mathbf{x}_{i,j+1}, \mathbf{x}_{i+1,j+1}, \mathbf{x}_{i,j+2}\right],\nonumber\\
S_5 &=& \left  [\mathbf{x}_{i-1,j-1}, \mathbf{x}_{i-2,j}, \mathbf{x}_{i-1,j}, \mathbf{x}_{i,j}, \mathbf{x}_{i-1,j+1}\right]. \nonumber
\eea
In this example, the total stencil $S$ is constructed to satisfy
$\displaystyle\bigcap_{m=1}^{5} S_m = \{\mathbf{x}_{i,j} \}$ 
and $\displaystyle\bigcup_{m=1}^{5} S_m = S$,
containing 13 data points
whose local indices range from $i-2,j-2$ to $i+2,j+2$,
excluding the 12 cells in the corner regions.
See Fig.~\ref{fig:2dS}, for a detailed schematic of the multi-substencil approach.  

Using these data points, we build the 13$\times$5 over-determined system
\begin{equation}
\begin{pmatrix} w_{1,1} 	& 0       	& 0       	 & 0         	& 0       \\
	         	w_{1,2} 	& w_{2,1} & 0       	 & 0        	 & 0       \\ 
	                w_{1,3} 	& 0       	& w_{3,1}  & 0        	 & 0       \\ 
	         	w_{1,4} 	& 0       	& 0        	& w_{4,1}   & 0       \\ 
			0       	& w_{2,2} & 0       	 & 0         	& 0       \\ 
	        		0       	& w_{2,3} & w_{3,2}  & 0        	 & 0       \\ 
			w_{1,5} 	& w_{2,4} & w_{3,3}  & w_{4,2}   & w_{5,1} \\ 
			0       	& 0       	& w_{3,4}  & w_{4,3}   & 0       \\
			0 		& 0	  	& 0        	& w_{4,4}   & 0       \\ 
			0		& w_{2,5} & 0 	   	  & 0 	 	 & w_{5,2} \\ 
			0      		 & 0       	& w_{3,5}  & 0        	 & w_{5,3} \\
			0      		 & 0    	 & 0        	& w_{4,5}   & w_{5,4} \\ 
			0      		 & 0      	 & 0        	& 0         	& w_{5,5} 
\end{pmatrix} 
\begin{pmatrix}\gamma_1 \\ \gamma_2 \\ \gamma_3 \\ \gamma_4 \\ \gamma_5 \end{pmatrix}
= \begin{pmatrix}w_{1} \\ w_{2} \\ w_{3} \\ w_{4} \\ w_{5} \\ w_{6} \\ w_{7} \\ w_{8} \\ w_{9} 
\\ w_{10} \\ w_{11}  \\ w_{12} \\ w_{13}  \end{pmatrix},
\label{eq:2Dgam}
\end{equation} 
which is solved using the QR factorization method for least squares. 

Notice that both the pointwise SE kernel 
and the integrated SE kernel
in Section \ref{sec:gp_ptwise_prol} and Section \ref{sec:gp_vol_prol}
are both isotropic kernels.  Hence, every 
$\mathbf{K}_{m,\sigma}$ and $\mathbf{C}_{m,\sigma}$ are identical over each substencil, illustrating that the 
WENO combination weights (i.e., $\mathbf{w}_m^T$)
and GP model weights (i.e., $\mathbf{w}_*^T$ and $\mathbf{T}_*^T$) only need to be computed and saved once per level,
and reused later.

The nonlinear weighting approach of GP-WENO designed in this probabilistic way has proven to be robust
and accurate in treating discontinuities~\cite{reyes_new_2016, reyes2019variable}. Regardless, the nature of its non-linearity requires
the calculation of nonlinear weights are to be taken place over the entire computational domain, 
consuming an extra computing time. In this regard, one can save the overall computation if the GP-WENO weighting
could only be performed when needed, i.e., near sharp gradients, identified by a shock-detector.
In our GP formulation, we already have a good candidate for a shock-detector, that is, the GP-based $\beta_m$.
To meet this, we slightly modify Eq.~\eqref{eq:eig} to 
introduce an optional switching parameter $\alpha$, defined by
%
\begin{equation} 
\alpha = \frac{\sum\limits_{i=1}^{2D+1} \displaystyle\frac{1}{\lambda_i}
(\mathbf{v}_i^T\mathbf{f})^2}
{\mathbb{E}_{arith}^2[\mathbf{f}] + \epsilon_2 }.
\label{eq:alph}
\end{equation}  
Here, the data array $\mathbf{f}$ includes the $2D+1$ data solely chosen from 
the center most substencil, e.g., $S_1$ in Fig.~\ref{fig:2dS}, 
$\mathbb{E}_{arith}^2$ is the squared arithmetic mean 
over the sampled coarse grid data points over $2D+1$ sized substencil centered
at the cell $(i,j,k)$, or $S_1$, 
that is,
\begin{equation}
\mathbf{E}^2_{arith}[\mathbf{f}] = \left(\frac{1}{2D+1}\displaystyle\sum_{\mathbf{x}\in S_1}f(\mathbf{x})\right)^2,
\end{equation}
and finally, $\epsilon_2$ is a safety
parameter in case the substencil data values are all zeros. 
Notice that this is just a scaled version of the $\beta_m$ in Eq.~\eqref{eq:eig}
for the central substencil $S_1$. Since the $\sigma$ GP model
 is built with smooth data in mind 
prescribed by the smooth SE kernel, 
this parameter will detect ``unlikeliness" in of the data $\bf{f}$ with respect
to the GP model. 
Note that the
critical value of $\alpha$, called $\alpha_c$, will be based on the kernel chosen.   
Without the normalization by the 
squared arithmetic  mean, this factor will vary based on mean value of the data.
In this regard, dividing by the average value of the data, $\bf{f}$, helps to
normalize the factor without changing the variability detection. 
From the statistical interpretation of the GP model $\mathbf{E}^2_{arith}[\mathbf{f}]$ may be viewed as a likelihood measure for a GP that assumes uncorrelated data (i.e., $\bfK_{ij}=\delta_{ij}$), so that $\alpha$ becomes normalized relative to the likelihood of another model.

We choose a critical value, $\alpha_c$ so that shocks, and high variability in $\bf{f}$ are detected when 
$\alpha > \alpha_c$; smooth and low variability when $\alpha \leq \alpha_c$. We heuristically set $\alpha_c=100$ in this 
strategy. 
Using this $\alpha$ parameter we have a switching mechanism between 
the more expensive nonlinear multi-substencil GP-WENO method in this section
and the linear single-stencil GP model in Sections \ref{sec:gp_ptwise_prol} and
\ref{sec:gp_vol_prol}.

Using the multi-substencil GP-WENO method, there are generally 
$2D+1$ dot products of the stencil size for each prolonged point, 
$\prod_d r_d$. In patch-based AMR, even though refined grids are 
localized around the regions containing shocks and  
turbulence, there are often areas of smooth flow in every patch. 
The use of the switch $\alpha$ allows us to reduce the computational 
complexity to one dot product of the stencil size for each coarse stencil that has smooth data, therefore 
reducing the cost to one dot product of the stencil size for each prolonged point. This method is extremely useful in 
3D and when the refinement ratio is greater than $\mathbf{r} = 2$.  

We conclude this section by making a remark 
on one significant feature of GP which we do not explore in our current study.
The multi-substencil GP-WENO methods, 
outlined in~\cite{reyes_new_2016,reyes2019variable}, in smooth flows can variably increase/decrease the order of accuracy.
However, in the application for AMR prolongation there may be large grids to be refined, 
so the increased computational cost can become undesirable. 
Note that the linear single model GP interpolation is still $\mathcal{O}(\Delta x^3)$, 
and serves as a high-order accurate prolongation that often matches the order of accuracy of the simulation. 
Reyes et al.~\cite{reyes_new_2016, reyes2019variable} discuss how to vary accuracy as a tunable 
parameter within the GP methodology. The studies show that the GP radius $R$ of 
the stencil dictates the order of accuracy. The method illustrated in 
this paper utilizes a GP radius $R =1$ and is $\mathcal{O}(\Delta x^3)$, 
however if one uses $R=2$ we can retrieve a method
that is $\mathcal{O}(\Delta x^5)$. 

\section{Implementation}
\label{sec:method}
The multi-substencil GP-WENO prolongation method is implemented within the AMReX framework. 
The AMReX framework utilizes a hybrid C++/Fortran library with many routines
to support the complex algorithmic nature of patch-based AMR and state-of-the-art high performance
computing. As an example, the object-oriented nature of C++ is fully utilized 
to furnish simple data and workflow.
In AMReX there is a virtual base class called \textit{Interpolator}. 
This class has many derivations including \textit{CellConservativeLinear}, an object for the functions 
related to a cell-based conservative linear interpolation. The methods presented in the current work 
reside in the \textit{CellGaussian} class within the AMReX framework. This class 
constructs a GP object, which contains the model weights for each of the 
$\prod_{d = 1}^{D} r_d$ new points per cell as member data. 
When a simulation is executed in parallelized format, each MPI rank has an Interpolator class, 
by which it helps to
avoid unnecessary communication. 
Computationally, the order of execution is as follows:
\begin{enumerate} 
\item The refinement ratio and $\Delta \mathbf{x}$ are passed on to the construction of the GP object.
\item Build GP covariance matrices for both interpolation $\bK, \bK_m$ and shock detection $\bK_{m,\sigma}$ ($\mathbf{C},\mathbf{C}_m$ and $\mathbf{C}_{m,\sigma}$ for cell-averaged data) using the SE kernel, Eq.~(\ref{eq:sqrexpcov}) (Eq.~(\ref{eq:intsq}) for the integrated kernel).
\begin{itemize}
    \item $\bK$ and $\bK_m$ are used for prolongation and should be used with the $\ell$ hyperparameter. This should be on the order of the size of the stencil to match our model assumption that the data varies smoothly over the stencil. We adopt $\ell = 12\cdot \Delta$, where $\Delta = \min(\Delta x_d)$
    \item $\bK_{m,\sigma}$ is used in the shock detection through the GP smoothness indicators and should take $\sigma$ as the length hyperparameter. $\sigma/\Delta \sim 1.5 - 3.5$ corresponding to the typical shock width in high-order Godunov method simulations. 
\end{itemize}
\item Calculate the GP weights $\mathbf{w}_*$ ($\mathbf{z}_*$) for all $\prod_{d=1}^{D} r_d$ prolonged points using Eq.~(\ref{eq:pinterp2}) (Eq.~(\ref{eq:vinterp2})). These weights are calculated only once for every possible refinement ratio and stored for use throughout the simulation. 
\item Compute the eigensystem of $\bK_{m,\sigma}$ as part of building the shock-capturing model. The eigenvectors $\mathbf{v}_i/\sqrt{\lambda_i}$ are stored for use in calculating $\beta_m$ and $\alpha$. 
\item Solve for $\boldsymbol{\gamma}$ for each prolonged point using the weights 
from $S_m$ and $S$. Just as for the GP weights $\mathbf{\gamma}$ is only calculated once and stored for each possible refinement ratio. 
\item For each coarse cell, the switch parameter $\alpha$ is calculated and compared to $\alpha_c$ (we choose $\alpha_c=100$).
\begin{itemize} 
    \item When $\alpha < \alpha_c$ the data is determined to be smooth enough not to need the full nonlinear GP-WENO prolongation. Instead the GP-weights over the total stencil $\mathbf{w}_*$ ($\mathbf{z}_*$) may be used without any weighting.
	\item In the case that $\alpha < \alpha_c$ the points are prolonged using the nonlinear multi-substencil GP-WENO model 
	(e.g., one of the methods in Sections~\ref{sec:gp_ptwise_prol} and \ref{sec:gp_vol_prol}, plus 
	the nonlinear controls in Section~\ref{sec:GP-WENO}).
\end{itemize}  
\end{enumerate}  

If needed, the parameter $\alpha_c$ can be tuned to a different value 
to alleviate the GP performance relating to sensitivity to shock-detection. 
By lowering $\alpha_c$ shocks will be detected more frequently, 
leading the overall computation to increase since  GP-WENO will be activated on
an increased number of cells.
%
In most practical applications such a tuning would be unnecessary considering that strong shocks are fairly localized, 
and in such regions $\alpha$ would retain a value much larger than $\alpha_c$ anyway.
Therefore, the condition $\alpha > \alpha_c$ for nonlinear GP-WENO would be met most likely 
over a wide span of possible values of $\alpha_c$ users might set.
Nonetheless, the localized nature of shocks allows the 
computationally efficient linear GP model to be used in simulations 
that do not require the frequent shock handling mechanism. 
In what follows, we set $\alpha_c = 100$ in the numerical test cases presented in Section~\ref{sec:testing}.

\begin{figure}
\begin{center}
\includegraphics[height=8cm, trim={2cm 3cm 2cm 3cm},clip]{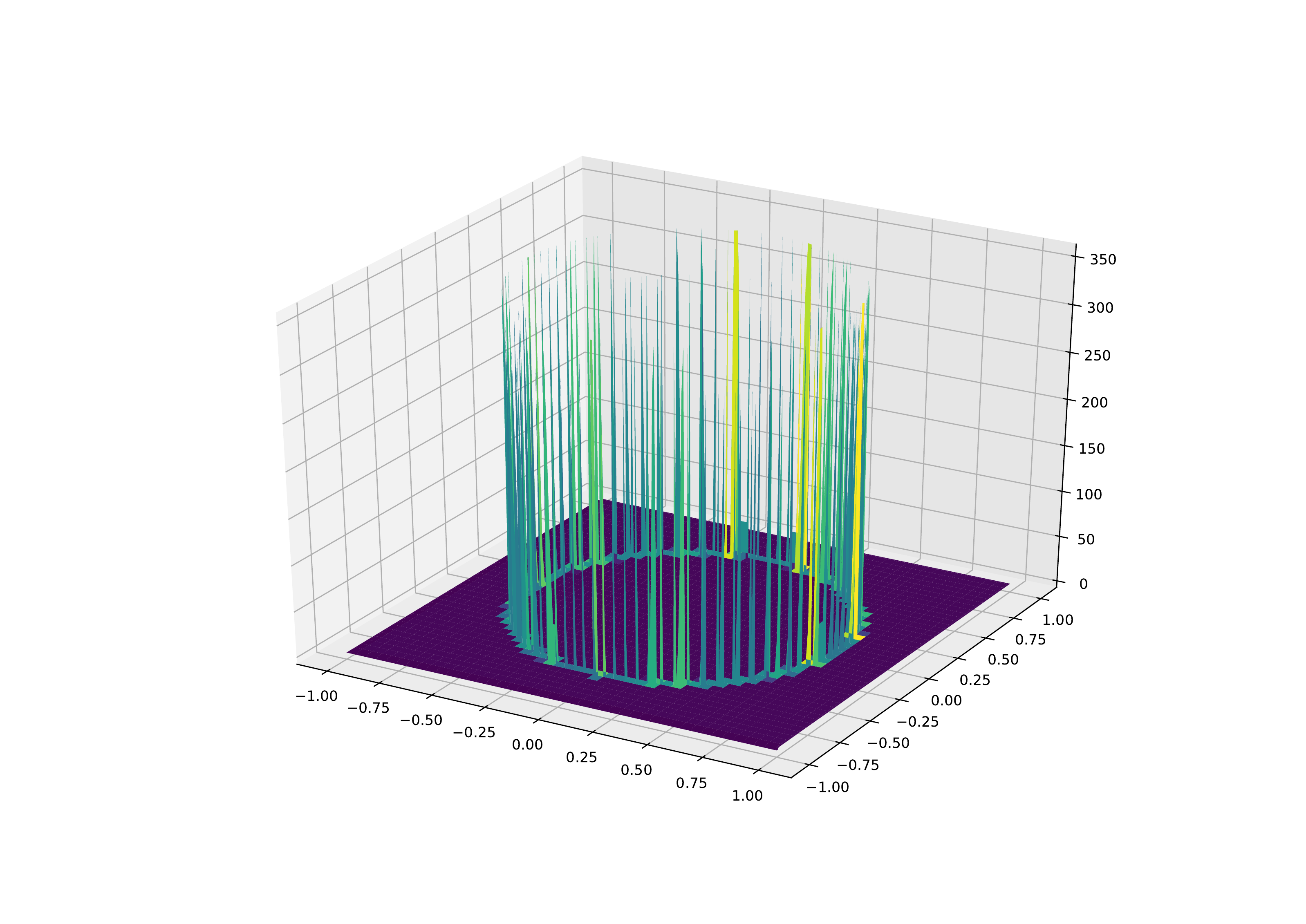}
\includegraphics[height=8cm, trim={2cm 3cm 2cm 3cm},clip]{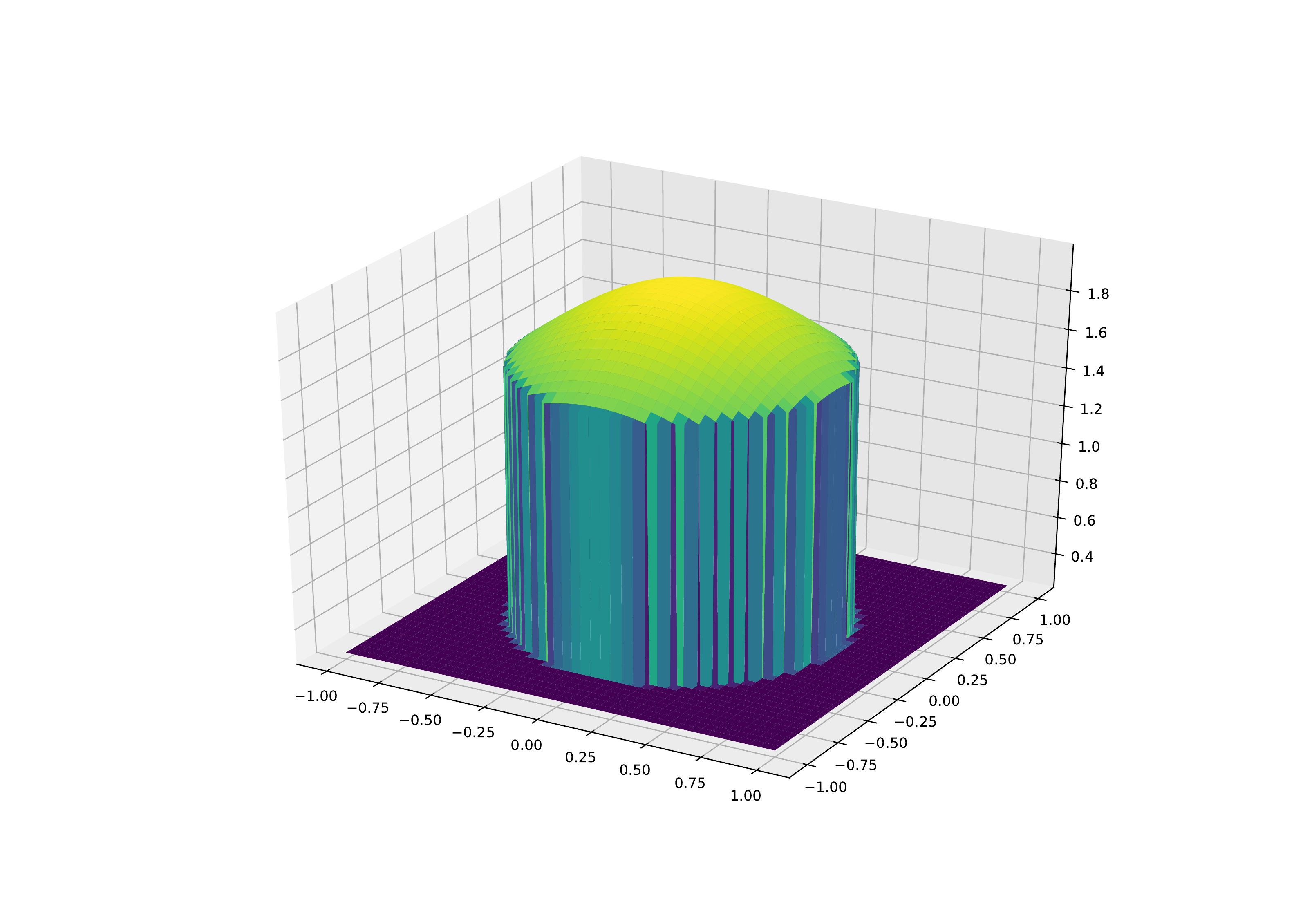}
\caption{\label{fig:alpha} The top plot displays the $\alpha$ values associated with the data from the 
function $f(x,y)$ depicted in the bottom plot.}
\end{center}
\end{figure}

To illustrate, 
we show the $\alpha$ values associated with a 
Gaussian profile elevated by the circular cylinder of height 0.25 defined by
the following function
\begin{equation}
f(x, y) = \begin{cases}
1 + \exp{\left(-(x^2 + y^2)\right) \quad \textrm{if}\hspace{1mm} (x^2 + y^2) < 0.5}, \\
0.25 \quad \textrm{else}.
\end{cases}
\end{equation}
In Fig.~\ref{fig:alpha},  we demonstrate how $\alpha$ varies over
the profile which combines the smooth continuous profile with the abrupt discontinuity.
It is observed that the $\alpha$ value is close to 2 over the continuous region. However, at the points 
corresponding to the sharp discontinuity, $(x^2 + y^2) = 0.5$, $\alpha$ soars to over 300,
resulting in the full engagement of the multi-substencil GP-WENO model near the discontinuity.
In the rest of the smooth region, $\alpha$ becomes much smaller, 
and therefore only using the linear GP model. 
This also tells us that the linear GP model would be a sufficient AMR prolongation algorithm
in an incompressible setting.

\subsection{AMReX Programming Directives} 

We make our implementation of the present work publicly available at \url{https://github.com/stevenireeves/amrex} in the GP-AMR branch. 
Written in C++, it utilizes AMReX's hardware-agnostic parallelization macros and lambda functions. The code is 
designed to utilize pragmas that declare the interpolation function as callable from either a CPU or GPU. The 
AMReX parallelization strategy is similar for both CPU-based supercomputers (e.g., Cori at NERSC) 
and GPU-based machines  (e.g., Summit at the Oak Ridge Leadership Co mputing Facility (OLCF), 
as well as Perlmutter at NERSC and the forthcoming Frontier at OLCF). 
The strategy is to use MPI for domain decomposition,
OpenMP for CPU based multi-threading, and CUDA (and HIP in the future) for GPU accelerators. 
Data allocation, CPU-GPU data transfers and handling are natively embedded 
in most AMReX data types and objects. 
For a more in-depth look into how the AMReX software framework is implemented we invite the interested readers 
to refer to \cite{amrex}. 

To provide a simple example into the AMReX style of accelerator programming, 
let us suppose that we wish to assign the integer value of 1 to the whole AMReX datatype 
\textit{Array4}. The datatype is a 
plain-old-data object which is a four dimensional array indexed 
as $( i, j, k, n)$. The first three indices are for spatial indices and 
the last one for each individual component (e.g., fluid density, $\rho$). 
We can use the AMReX lambda function \textbf{AMREX\_PARALLEL\_FOR\_4D} to expand a 4D loop in a parallel fashion.
For instance, \textbf{AMREX\_PARALLEL\_FOR\_4D} distributes the code segment in Listing 1 to an
equivalent format in Listing 2: 

\begin{center}
\begin{lstlisting}[frame=single, language=C++, basicstyle=\scriptsize, 
caption={AMReX Directive for Kernel Launch}]
AMREX_PARALLEL_FOR_4D(bx, ncomp, i, j, k, n, {
	my_array(i,j,k,n) += 1; 
}); 
\end{lstlisting} 
\end{center}


\begin{center}
\begin{minipage}{\linewidth}
\begin{lstlisting}[frame=single, language=C++, basicstyle=\scriptsize, 
caption={Expanded For Loop}]
for(int i = lo.x; i < hi.x; ++i){
	for(int j = lo.y; j < hi.y; ++j){
		for(int k = lo.z; k < hi.z; ++k){
			AMREX_PRAGMA_SIMD 
			for(int n = 0; n < ncomp; ++n) 
				my_array(i,j,k,n) += 1; 
		}
	}
}
\end{lstlisting}
\end{minipage}
\end{center}

This formulation allows for one code to be compiled for either CPU running or GPU launching. 
The AMReX lambda functions are expanded by the compiler, and the box dimensions (lo.x - hi.x, etc) 
are different based on the target device. For GPUs the lo and hi variables are set based on how much 
data each GPU thread will handle. In the CPU version, the lo to hi are the dimensions of target tile boxes respectively. 
Essentially the lambda handles the GPU kernel launch or CPU for-loop expansion for the developer/user.

Assuredly, there are other approaches available to launch a parallel region in AMReX for GPU extension.
We recommend that further interested reader view the GPU tutorials in AMReX source code 
for more information on various types of launch macros~\cite{amrex}.

\begin{figure}
\begin{center}	
\includegraphics[width=0.65\textwidth]{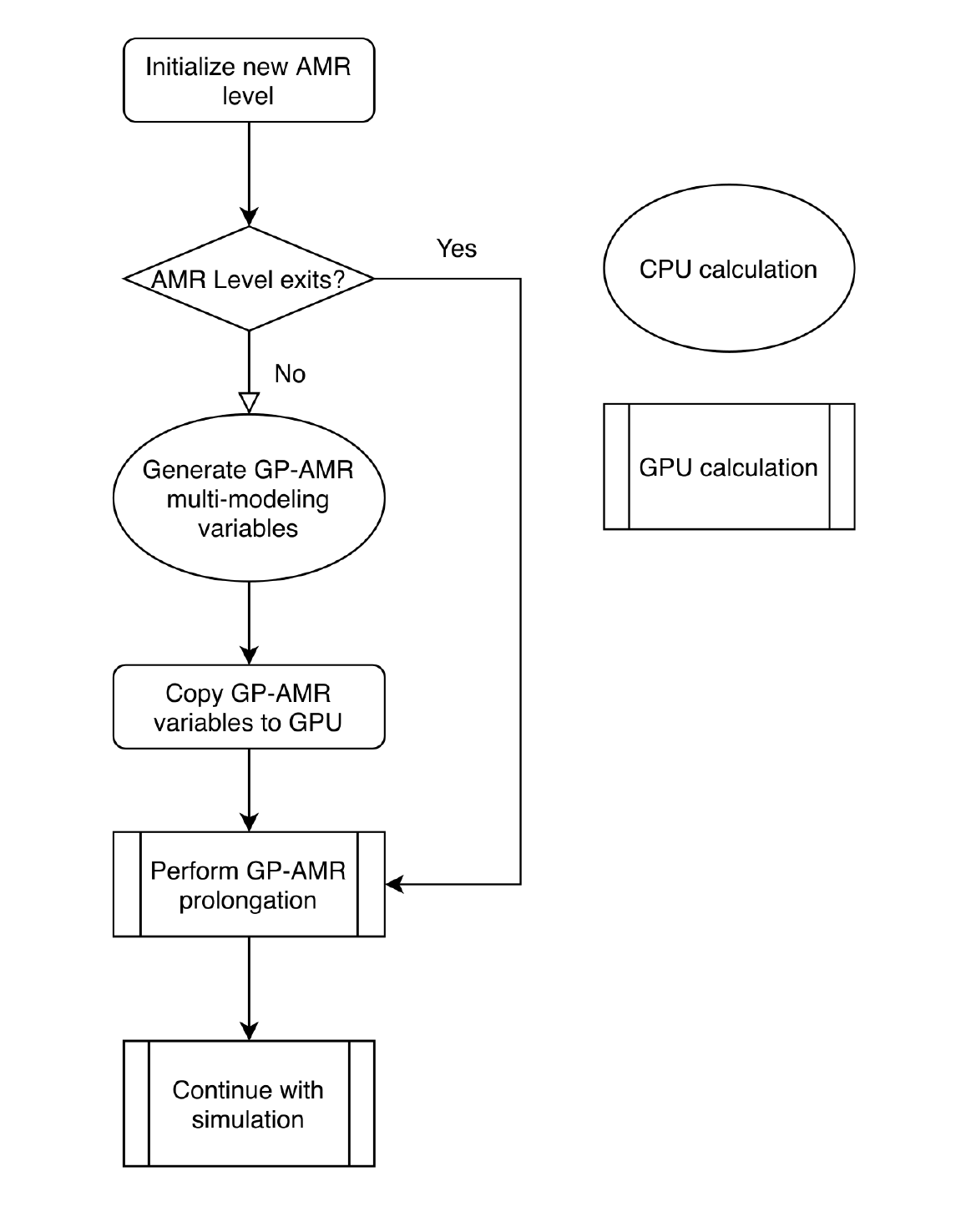}
\caption{\label{fig:diagram} Diagram illustrating a call graph for GP-AMR utilizing GPUs as accelerators.
} 
\end{center}
\end{figure}

While a detailed description of GPU computing is beyond the scope of this paper, 
we provide a general principle of our strategy to implement the GP-AMR algorithm for GPUs: 
\begin{enumerate}
\item Construct the model weights $\mathbf{T}^T_{*}\mathbf{C}^{-1}$ for each stencil $S_m$, $\boldsymbol{\gamma}$ and the 
eigensystem of $\mathbf{C}_\sigma$ on the CPU at the beginning of program execution or at the initialization of each AMR level.
\item Create a GPU copy for these variables and transfer them to the GPU global memory space. 
Every core on the GPU will need to access them, but do not need their own copy. 
\item Create a function for the prolongation. This function will require both the coarse grid data 
as an input, and the fine grid data as an output. Both arrays will need to be on the global GPU memory space.
This function will be launched on the GPU, and the fine level will be filled accordingly.  
\end{enumerate}
In general, with GPU computing, it is best to do as few memory transfers between the CPU and GPU as possible
because a memory transfer can cost hundreds or thousands of compute cycles and can drastically slow down an application. 
To further explain these steps, Figure~\ref{fig:diagram} is of an example call graph along with CPU-GPU memory transfers. 
In this diagram, it is already assumed that the course and fine state variables have been constructed and allocated on the GPU 
respectively, as is with the case in AMReX.  

\section{Results}
\label{sec:results}
In this section, we present the performance of the new GP-based prolongation model compared with the 
default conservative linear polynomial scheme in AMReX.
To illustrate the utility of the new GP-based prolongation scheme in fluid dynamics simulations,
we integrated the prolongation method in two different AMReX application codes,
Castro~\cite{castro} -- a massively
parallel, AMR, compressible Astrophysics simulation code,
PeleC -- a compressible combustion code~\cite{pelec}, as well as 
a simple advection tutorial code built in AMReX. 

\subsection{Accuracy}
To test the order of accuracy of the proposed method, 
a simple Gaussian profile 
is refined with the GP prolongation method. This profile follows the formula
\begin{equation}
f(\mathbf{x}) = \exp(-||\bx||^2)
\label{eq:acc}
\end{equation}
where $\mathbf{x} \in [-2, 2] \times [-2, 2]$. 
We compare the prolonged solution, denoted as $f_p(\bx)$, against the analytical value, $f(\bx)$,
associated with the Gaussian profile function. 
We find that the accuracy of the cell-averaged GP prolongation routine matches with 
the analysis in~\cite{reyes_new_2016,reyes2019variable}.
The convergence rate of the error in 1-norm, $E=||f-f_p||$, 
computed using the GP prolongation model with $R=1$,
exhibits the expected third-order accuracy, following the theoretical slope of third-order convergence
on the grid scales, $\mathcal{O}(\Delta x^3)$. 

\begin{figure}
\begin{center}
\includegraphics[width=0.65\textwidth]{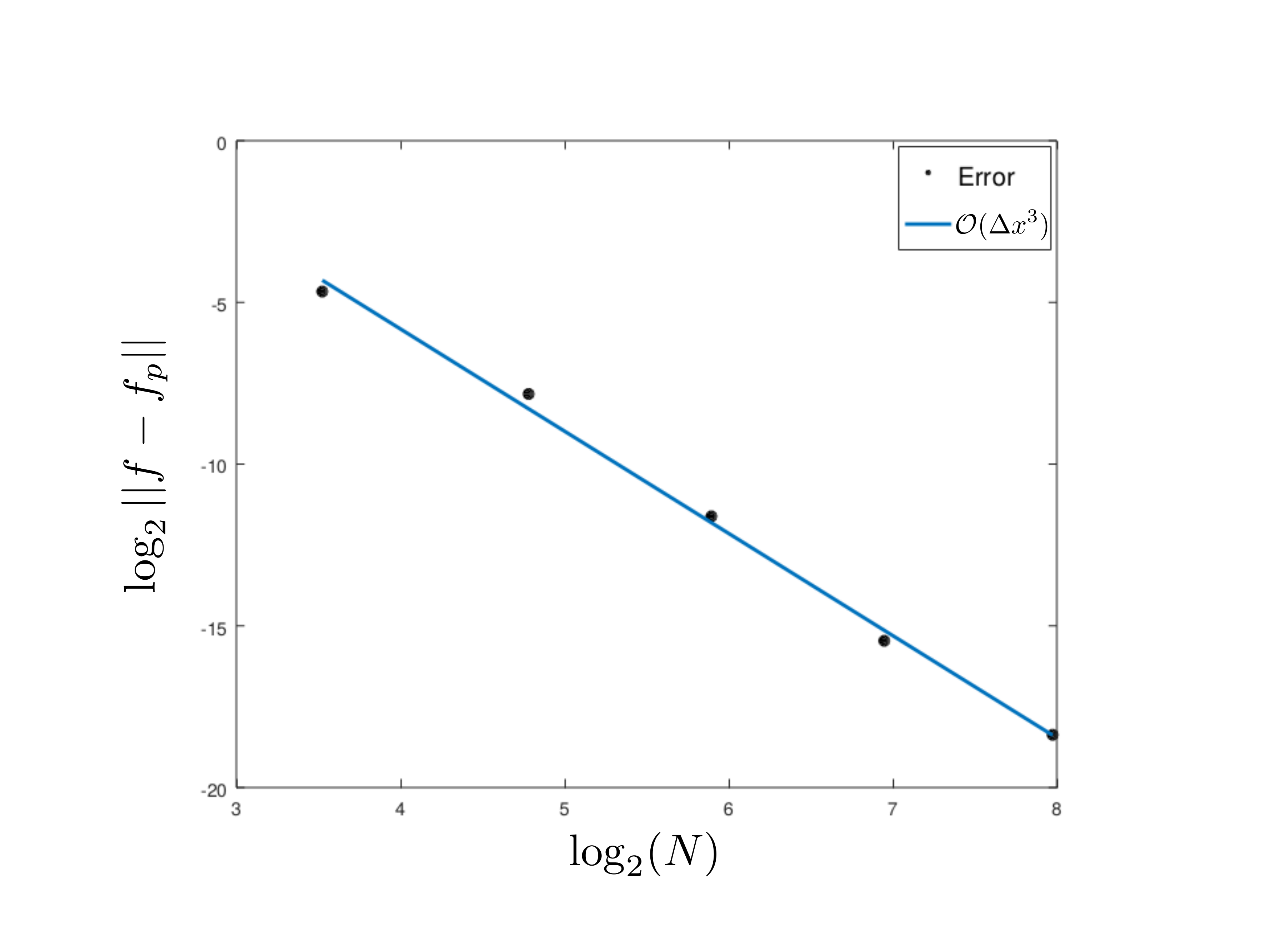}
\caption{\label{fig:err} Convergence for the GP method. 
The quantities are measured in $\log$ of base 2 to better cope with the refinement jump ratio of 2.}
\end{center}
\end{figure}

\subsection{GP-AMR Tests} 
\label{sec:testing}
There are several test problems we will examine. 
First will be a single vortex advection, provided as the \textit{Advection\_AmrLevel} tutorial in 
AMReX. Next we will present a modified version of the slotted cylinder problem from~\cite{mapped}.  
 The subsequent problems using Castro are some classic hydrodynamic test problems including 
the Sedov implosion test~\cite{sedov} and
the double Mach reflection problem~\cite{dmr}.
Lastly, we discuss a premixed flame simulation from PeleC.
In all tests we use $\ell = 12\cdot\min(\Delta x_d)$. For the 2D simulations, $\sigma = 3\cdot\min(\Delta x_d)$ is used, and for the 
3D simulations we use $\sigma = 1.5\cdot\min(\Delta x_d)$. 
\subsubsection{Single Vortex Advection using the AMReX tutorial} 
The first test is a simple reversible vortex advection run. A radial profile 
is morphed into a vortex and reversed back into its original shape. 
This stresses the AMR prolongation's ability to recover the  profile after it has been advected into the coarse cells so that at the final time the solution can return to it's original shape. 
The radial profile initially is defined by
\begin{equation}
f(x, y) = 1 + \exp\Bigl[-100\left(\left(x-0.5\right)^2 + \left(y-0.75\right)^2\right)\Bigr].
\end{equation}
The profile is advected with the following velocity field:
\begin{equation}
\mathbf{v}(x,y,t) = \nabla \times \psi 
\label{eq:veloc}
\end{equation}
which is the curl of the stream function 
\begin{equation}
\psi(x,y,t) = \frac{1}{\pi}\sin\left(\pi x\right)^2 \sin\left(\pi y\right)^2 \cos\left(\pi \frac{t}{2}\right) 
\end{equation}
Here $(x,y)\in[0,1]\times[0,1]$. 
In this demonstration, the level 0 grid size is $64 \times 64$, and has
two additional levels of refinement surrounding the radial profile. The simulation is an 
incompressible advection problem using the Mac-Projection to compute the incompressibility condition
enforcing the divergence-free velocity fields numerically, $\nabla \cdot \mathbf{v}=0$~\cite{mac}. 
The flux is calculated by a simple second-order accurate upwind linear reconstruction method.
Although the overall solution is second-order which is lower than the third-order
accuracy of the GP prolongation method,
this example still illustrates the computational performance of the GP-based method over 
the default conservative second-order prolongator.  

\begin{figure}
\begin{center}
\begin{subfigure}{0.33\textwidth}
\includegraphics[width=\textwidth]{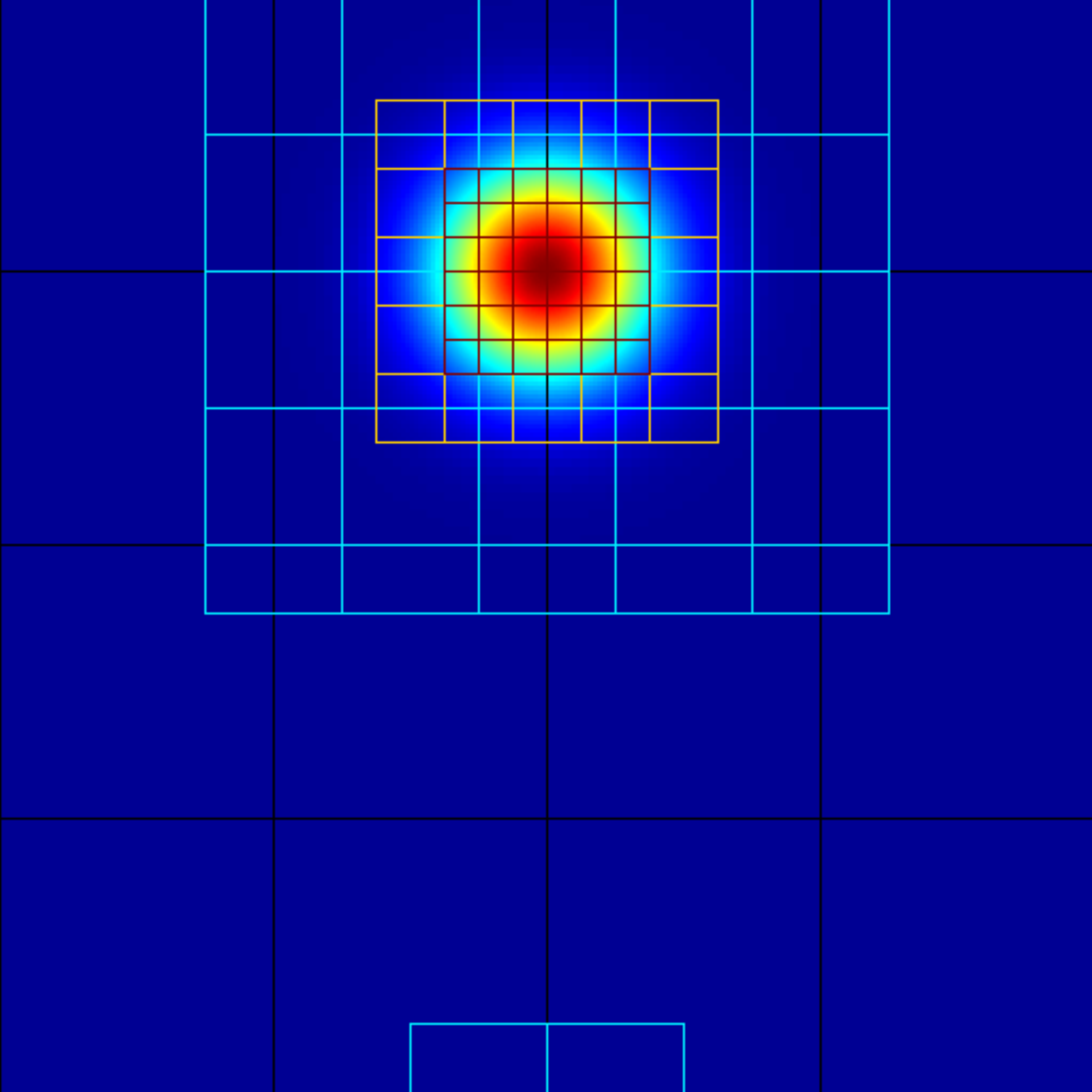}
\caption{$t = 0$}
\end{subfigure}%
\hfill
\begin{subfigure}{0.33\textwidth}
\includegraphics[width=\textwidth]{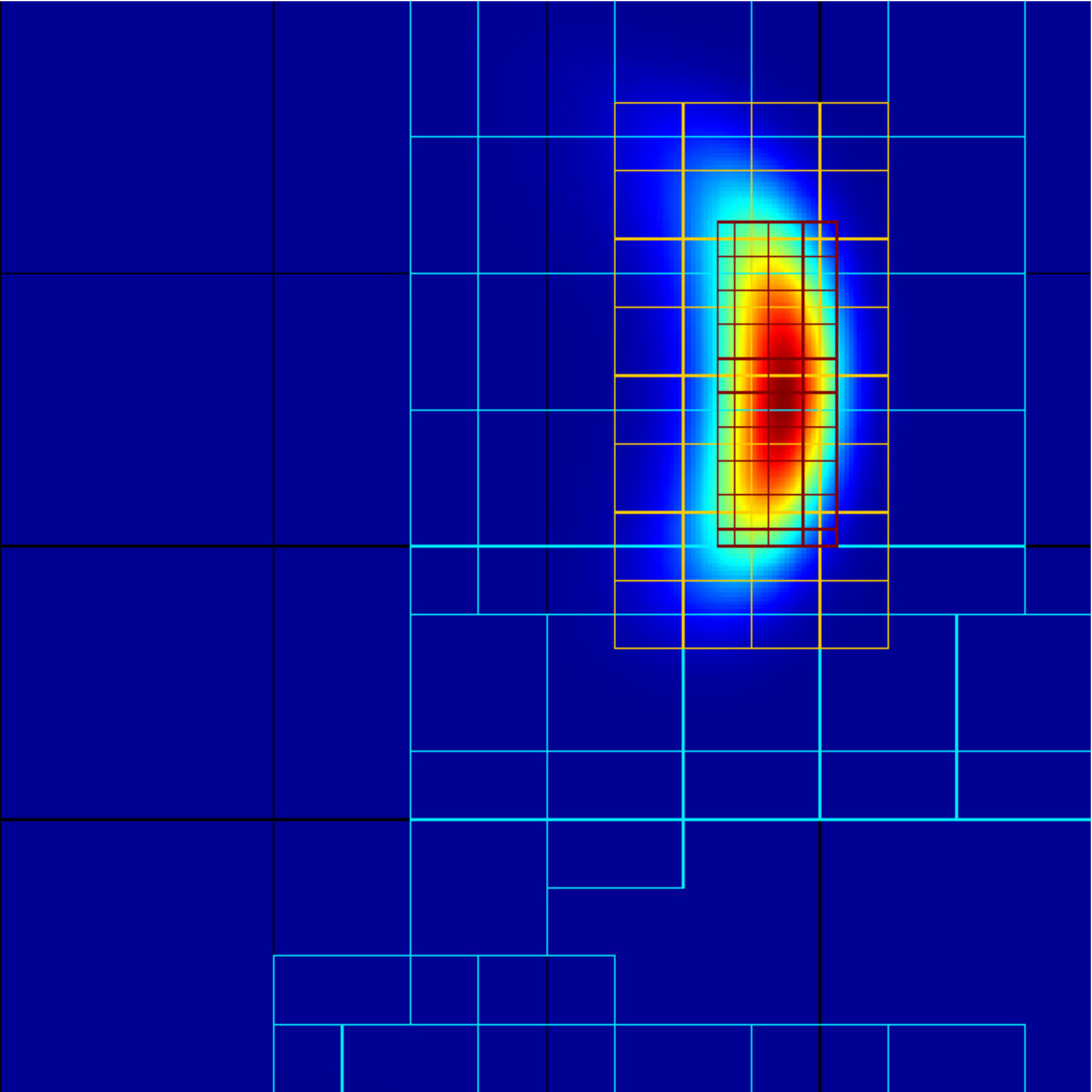}
\caption{$t = 0.282$}
\end{subfigure}%
\hfill
\begin{subfigure}{0.33\textwidth}
\includegraphics[width=\textwidth]{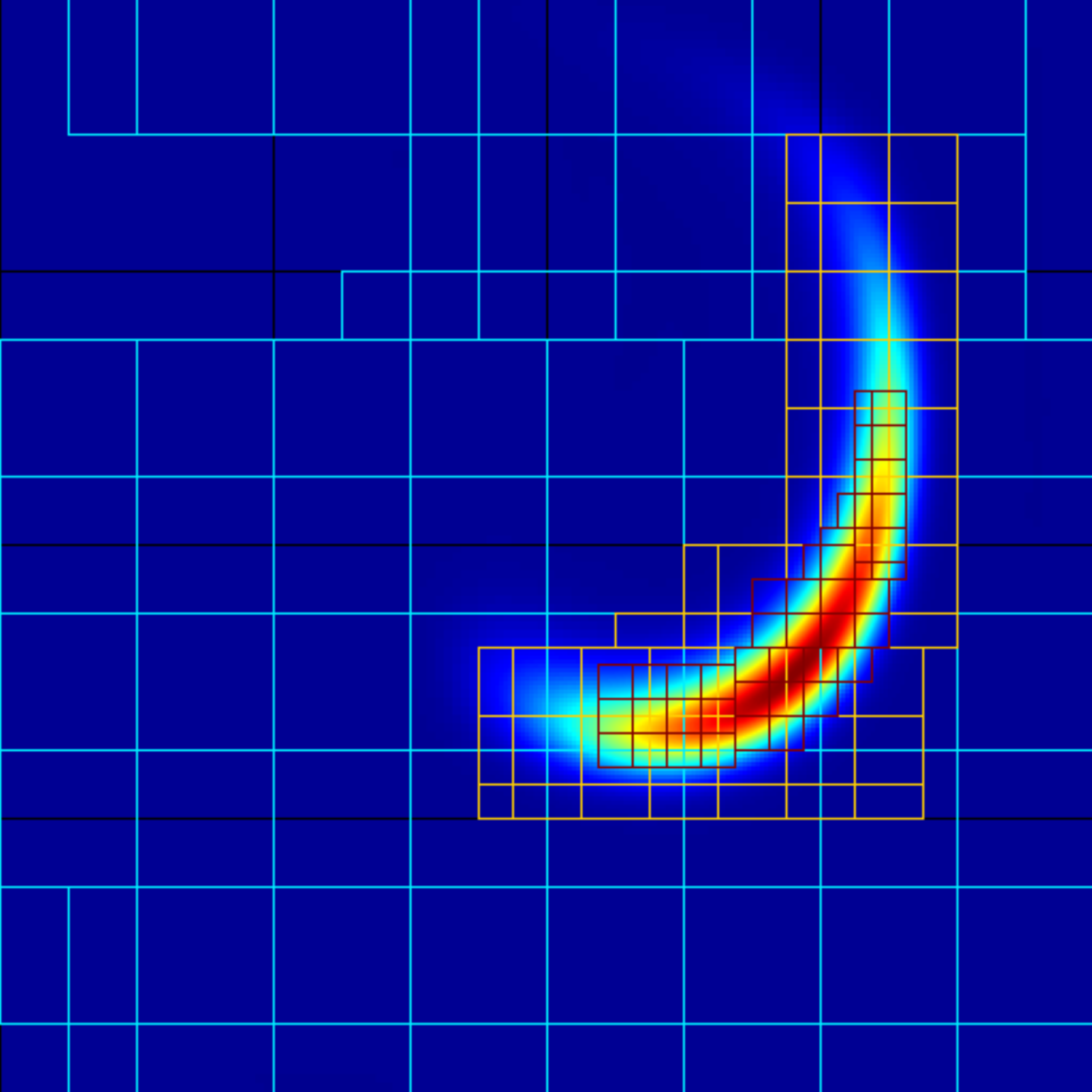}
\caption{$t = 0.651$}
\end{subfigure}

\bigskip
\bigskip

\begin{subfigure}{0.33\textwidth}
\includegraphics[width=\textwidth]{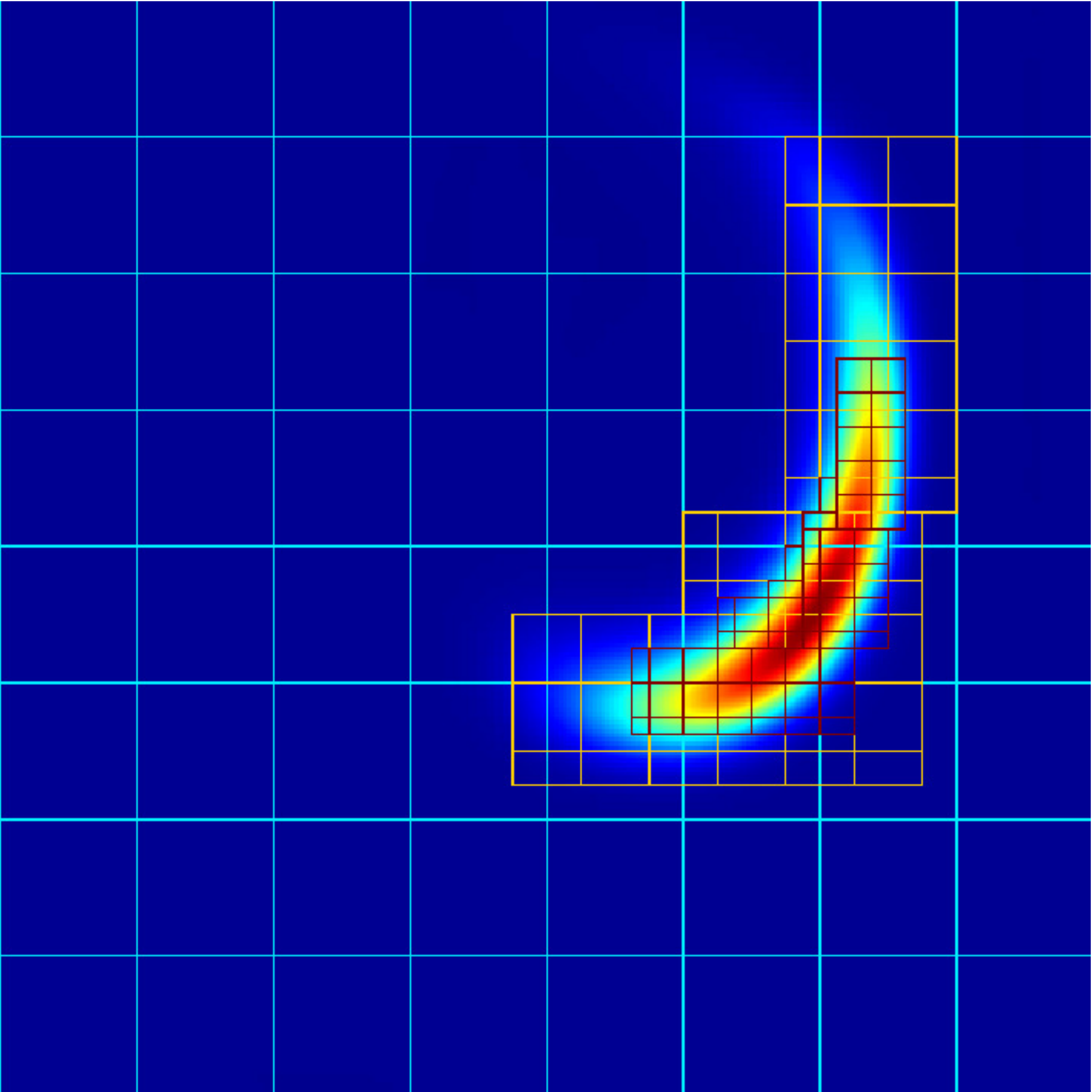}
\caption{$t = 1.447$}
\end{subfigure}%
\hfill
\begin{subfigure}{0.33\textwidth}
\includegraphics[width=\textwidth]{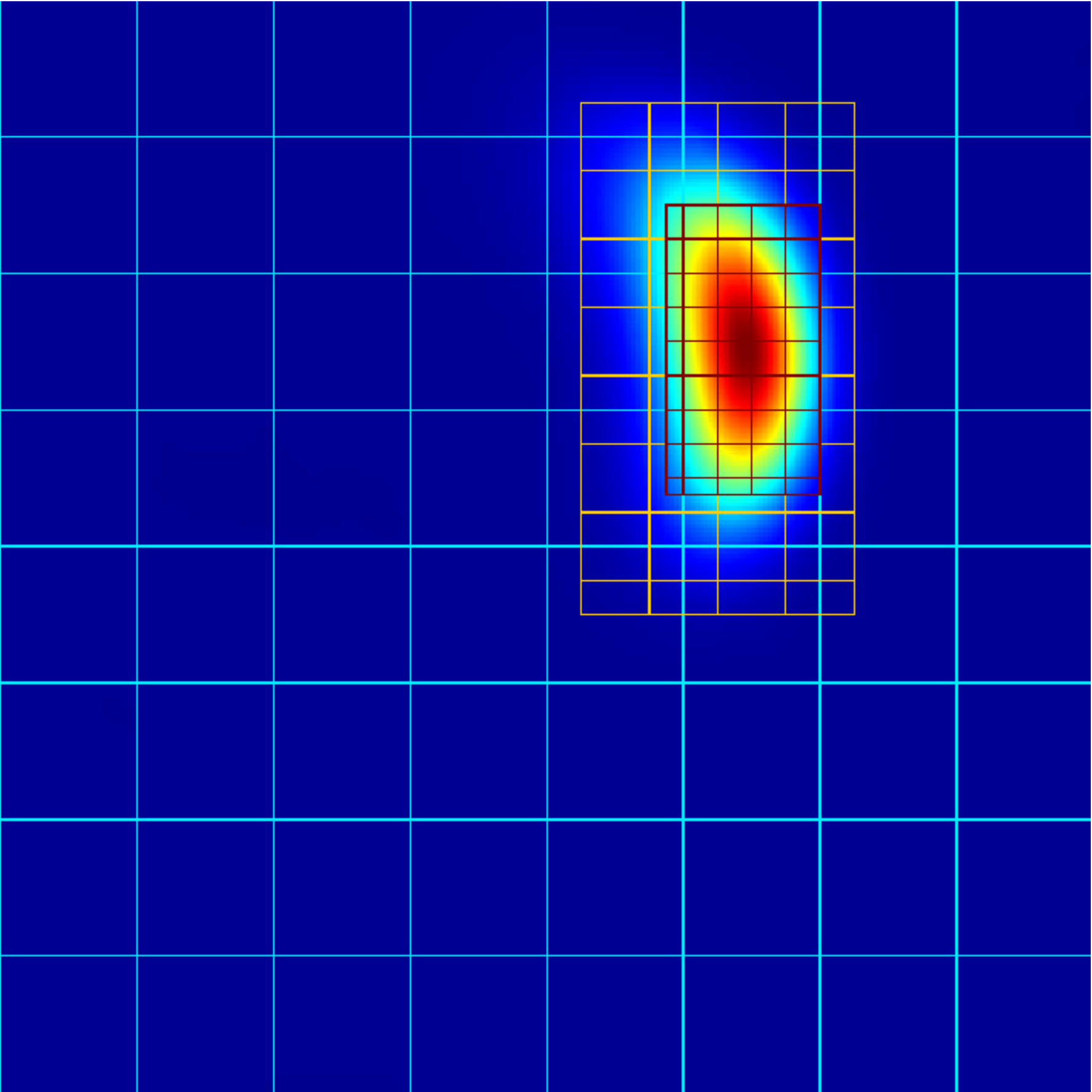}
\caption{$t = 1.785$}
\end{subfigure}%
\hfill
\begin{subfigure}{0.33\textwidth}
\includegraphics[width=\textwidth]{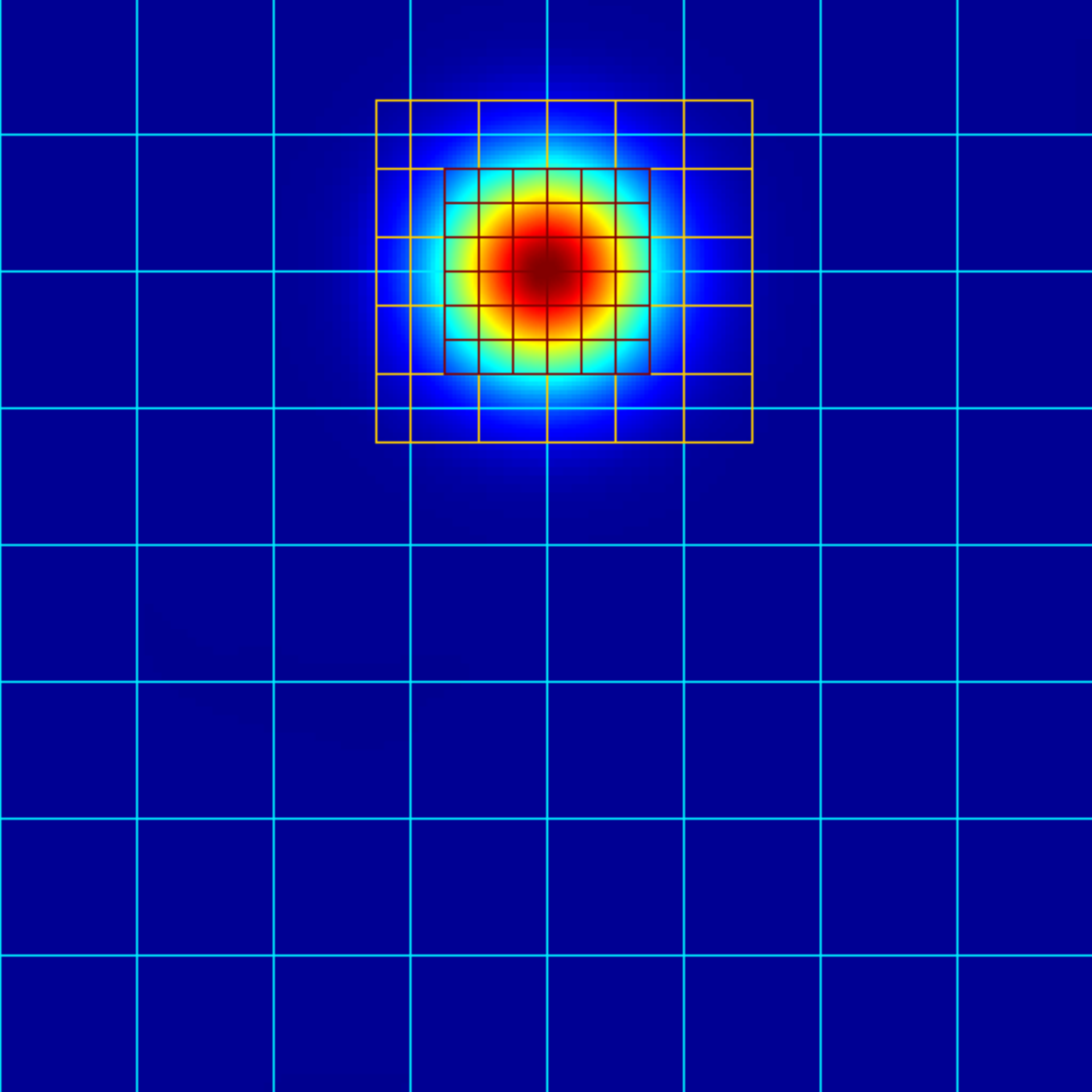}
\caption{$t = 2$}
\end{subfigure}
\caption{\label{fig:amrlev} The progression of the 2D radial profile with 4 
levels of refinement using the multi-substencil GP prolongation algorithm.
} 
\end{center}
\end{figure} 

The simulation is finished at $t = 2$.
We used sub-cycling of time-steps to improve the overall performance, 
in which a smaller time-step $\Delta t_f$ is used on a finer level 
to advance the regional solutions for stability. 
The coarser level solutions which advance with a larger time-step $\Delta t_c$ 
await until the solutions on the finer
levels catch up with the global simulation time $t_g=t^{n}+\Delta t_c$ over the number of sub-cycling
steps $N_{\mbox{subcycle}}=\Delta t_c/ \Delta t_f$.
We present the performance and accuracy results for this problem in Table~\ref{tab:SingleVort}.
 \begin{table}
  \begin{center}
    \caption{Accuracy and Performance of GP-AMR against the default linear AMR for the 
	single vortex test 
	     on a workstation with an Intel i7-8700K processor, with 6 MPI ranks.}
    \label{tab:SingleVort}
    \begin{tabular}{l|c|c|c|c} 
       & Execution Time & Prolongation Time & \# of calls  & $L_1$ error \\
	\hline
      $2D$ GP-AMR & 0.2323s & 0.004168s & 9115 & 0.00033 \\
      $2D$ Linear & 0.2335s & 0.008436s & 9113 & 0.00071 \\
      $3D$ GP-AMR & 1.6523s & 0.086361s & 21929 & 0.00151 \\
      $3D$ Linear & 1.6640s & 0.157623s & 21893 & 0.00160 \\
     \end{tabular}
  \end{center}
\end{table}

Since the two methods are of different order, 
they can yield different AMR level patterns
which can lead to the slight difference in the number of function calls. 
In the prolongation functions we find that the default linear prolongation took approximately 
twice as much time than the GP prolongation. This is due to the smoothness of the solution not requiring 
the multi-modeled treatment, allowing for the simplified GP algorithm to be used.
However, the overall simulation times  
were equally comparable,
since there were larger areas (or more cells) that followed the profile
and were computed in the finest AMR level in the GP case than the linear case.
We note that the cost of computing the GP model weights were negligible in comparison 
to the program's execution time of 
0.0002306 seconds on average, being called twice (since 
there were two levels) per MPI rank. We also check the level averaged $L_1$ error between the solution at 
$t=2$ and the solution at $t=0$ for both AMR prolongation methods. In the 2D case we find that the GP-AMR 
solution is approximately half of the error produced by the default method. This highlights the utility of a high-order prolongation method, as smooth features are better recovered after they have been advected into coarser cells.

Another useful examination is the analogous problem in 3D in which the computational stencils 
for both methods grow. For the 3D version we use a 
$32 \times 32 \times 32$ base grid with 2 levels of refinement. The details of this simulation can also 
be found in Table~\ref{tab:SingleVort}. We note that
a parallel copy operation becomes slightly more expensive with GP
because the need for the GP multi-substencil grows
on non-smooth regions to handle discontinuities in a stable manner,
as managed by the $\alpha_c$ parameter. 
This becomes more apparent in 3D, as the computational stencil effectively grows from 
7 to 25 cells when using the multi-substencil approach. In this 3D benchmark, the difference in error 
between these methods is less than in the 2D case. The GP-AMR simulation still outperforms the linear
prolonged simulation, but to a smaller degree. The 3D simulation has a coarser base grid 
by a factor of 2 than the 2D simulation, leading to a coarser finest grid in the simulation.
This is why the $L_1$ errors are greater in the 3D case. 

In Figure~\ref{fig:amrlev} we show the same 2D single vortex advection with GP-AMR
on a base grid of $64 \times 64$ with 4 levels of refinement
to illustrate the GP method with a more production level grid configuration.  

\subsubsection{Slotted Cylinder as another AMReX Test}
Another useful test is the slotted cylinder advection presented in~\cite{mapped}. 
In this paper we do not use an exact replica of this problem, but instead we utilize the initial profile
and perform a similar transformation as in the previous problem. That is, the slotted cylinder is morphed
using the same velocity used in the Single Vortex test.
Because the advected profile is now piecewise constant, compared to the smooth profile in the previous problem, it will require the use of the nonlinear GP-WENO prolongation and test the capabilities of the GP-based smoothness indicators for prolongation. 
 
\begin{figure}
\begin{center} 
\includegraphics[width=8cm]{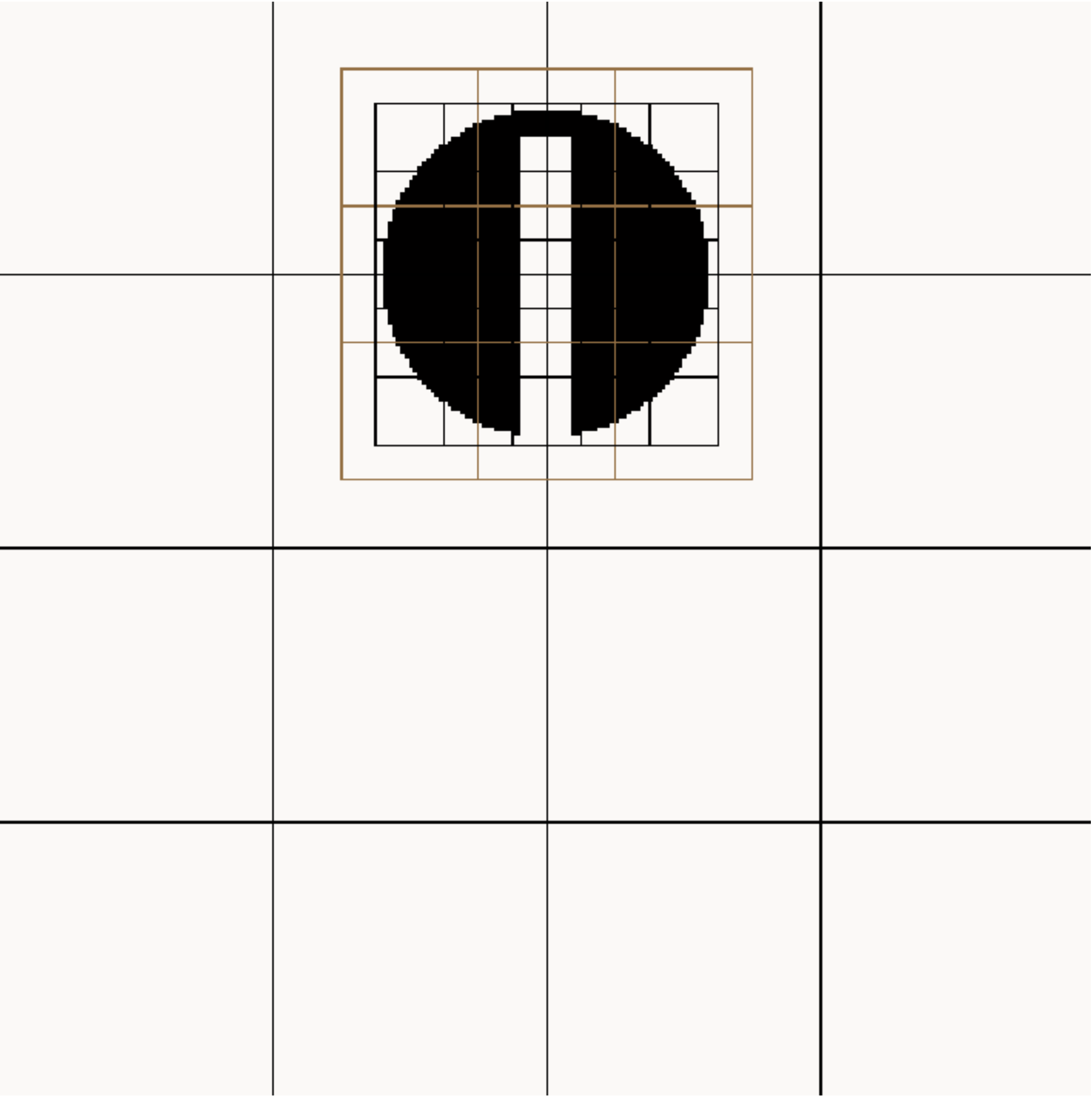}
\caption{\label{fig:slot0} The slotted cylinder at $t=0$ over the entire domain with 3 AMR levels. }
\end{center}
\end{figure}

The slotted cylinder is defined as a circle (in 2D) of radius $R=0.15$ centered at $\bx_c = (x_c, y_c) = (0.5, 0.75)$ with a slot of width $W=0.05$ and height $H=0.25$ removed from the center of the cylinder. 
The initial condition is given by
\[\phi_0(\mathbf{x}) =  \begin{cases}
0, \quad \mbox{if } R < \sqrt{(x-x_c)^2 + (y-y_c)^2}, \\
0, \quad \mbox{if } |2x_c| < W \;\; \mbox{and } 0 < y_c + R < H, \\
1, \quad \textrm{else},
\end{cases}\]
where $(x,y)\in[0,1]\times[0,1]$. The initial profile is shown in Figure~\ref{fig:slot0}.

\begin{figure}
\begin{center}
\begin{subfigure}{\textwidth}
\includegraphics[width=0.325\textwidth]{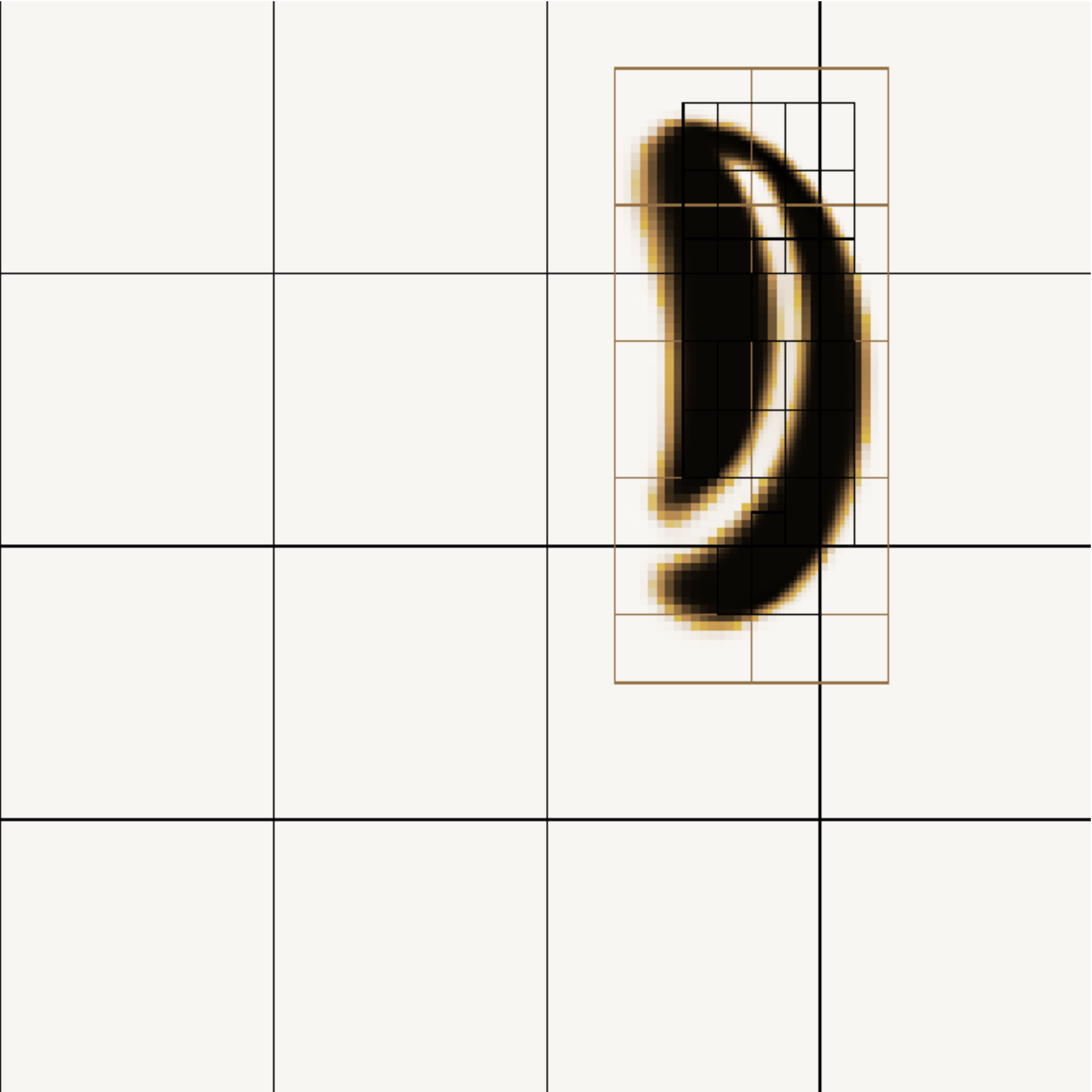}
\hfill
\includegraphics[width=0.325\textwidth]{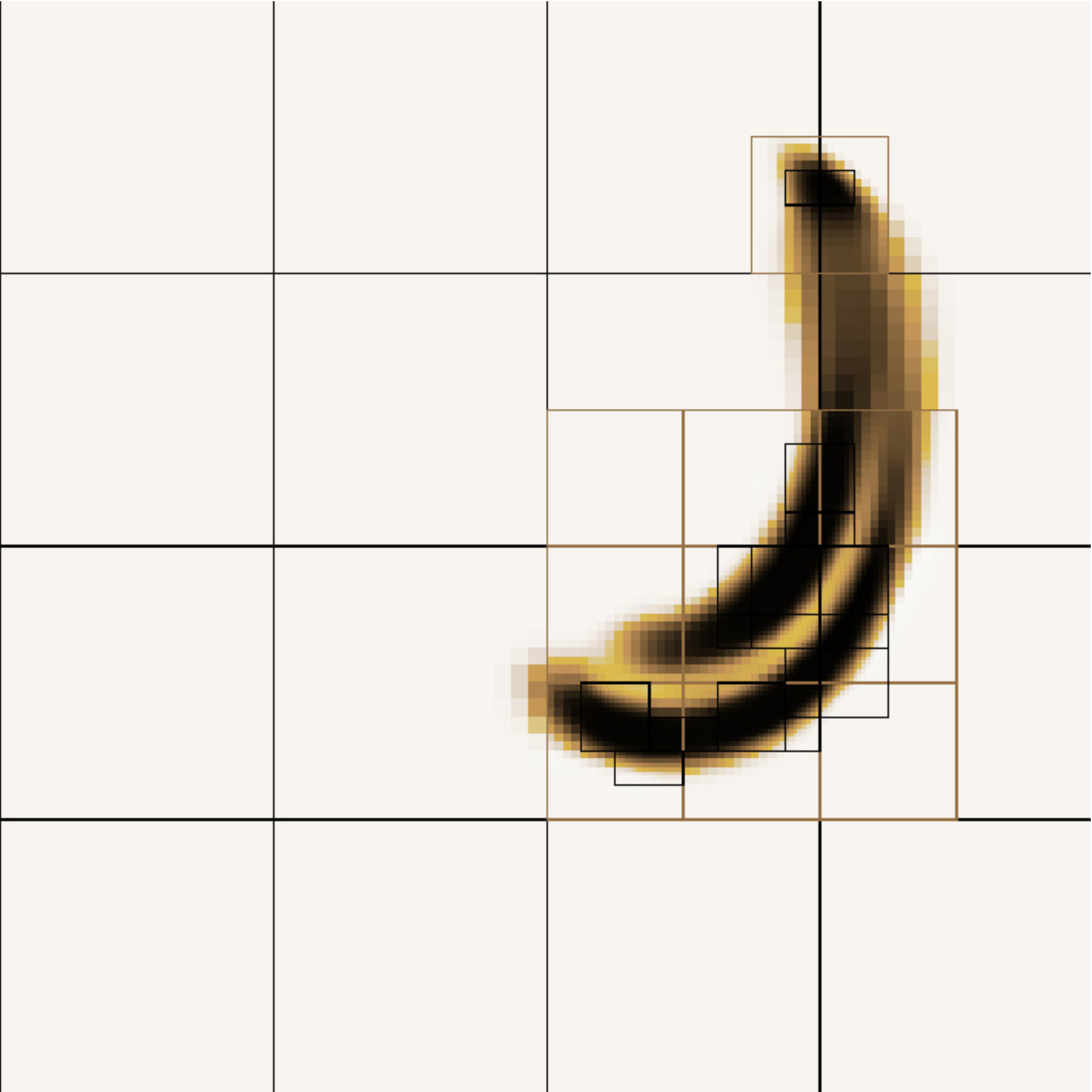}
\hfill
\includegraphics[width=0.325\textwidth]{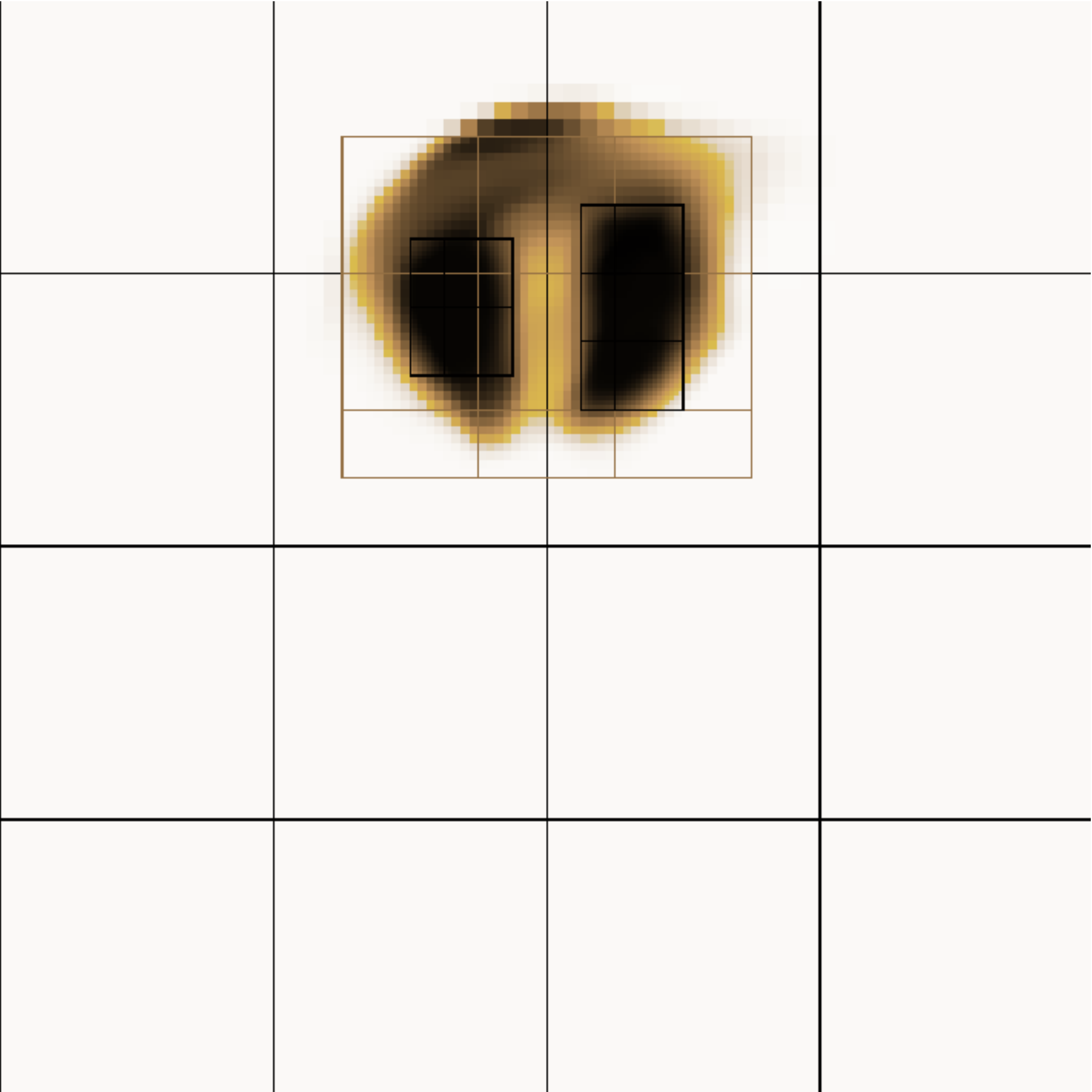}
\caption{Linear AMR prolongation}
\end{subfigure}
\\
\bigskip\bigskip
\begin{subfigure}{\textwidth}
\includegraphics[width=0.325\textwidth]{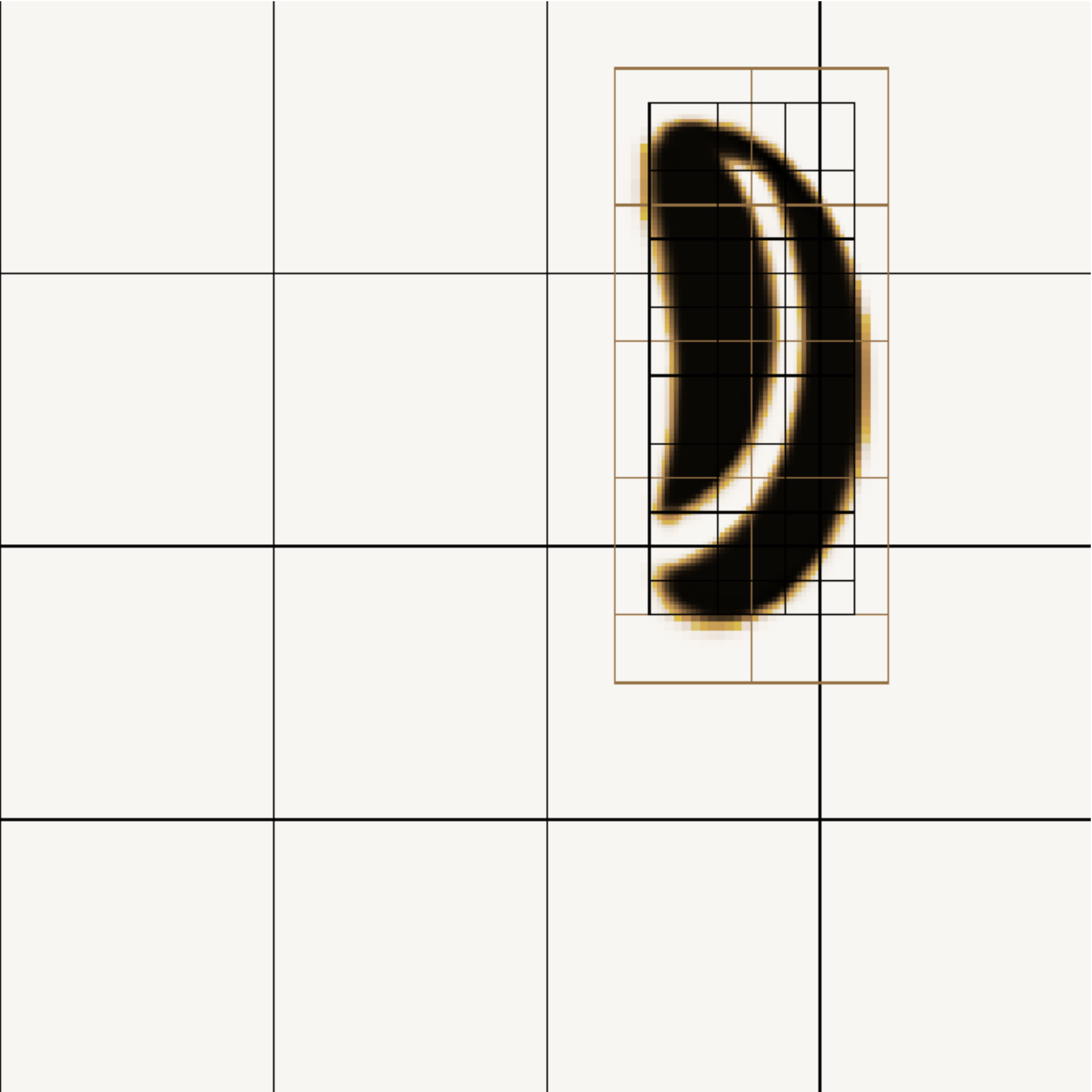}
\hfill
\includegraphics[width=0.325\textwidth]{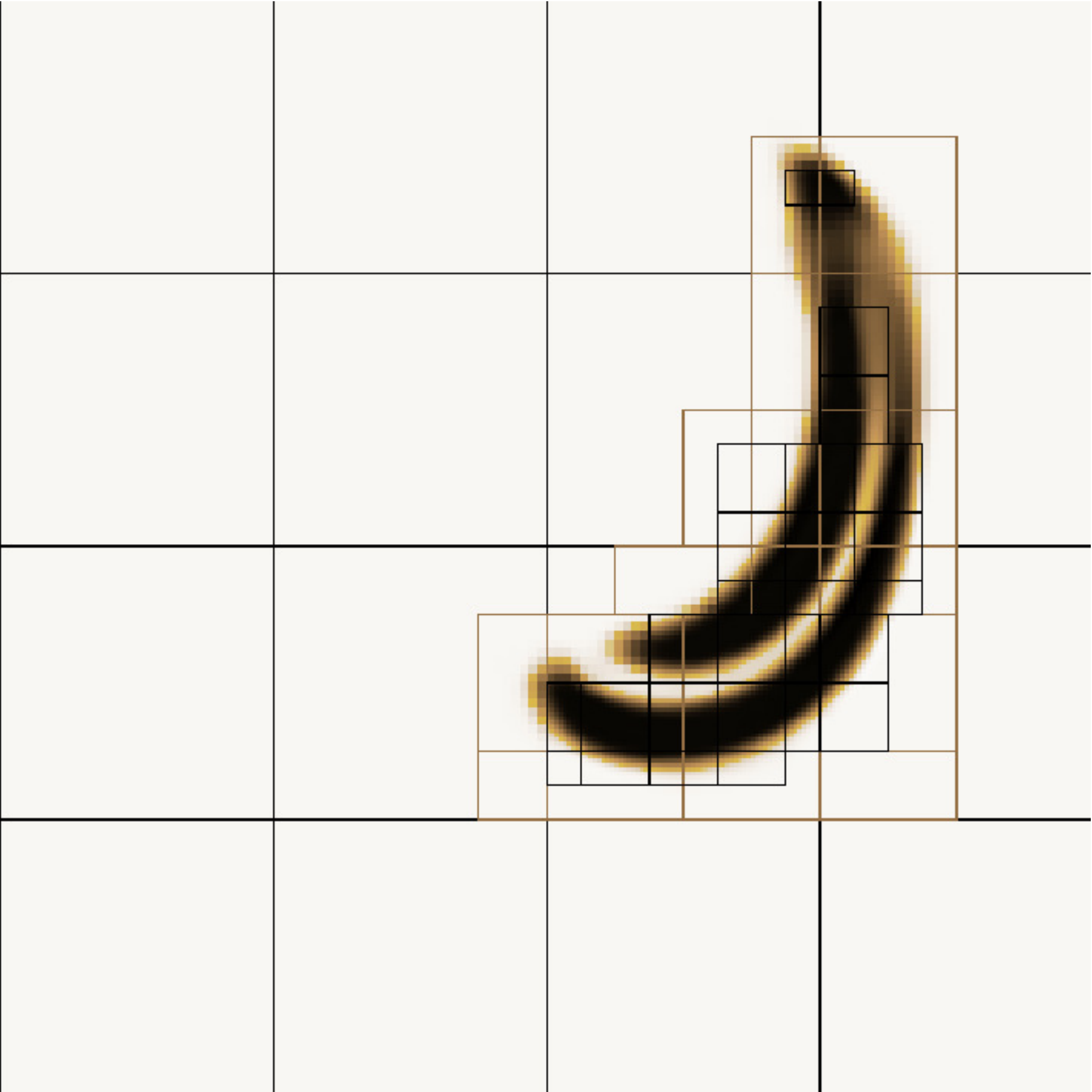}
\hfill
\includegraphics[width=0.325\textwidth]{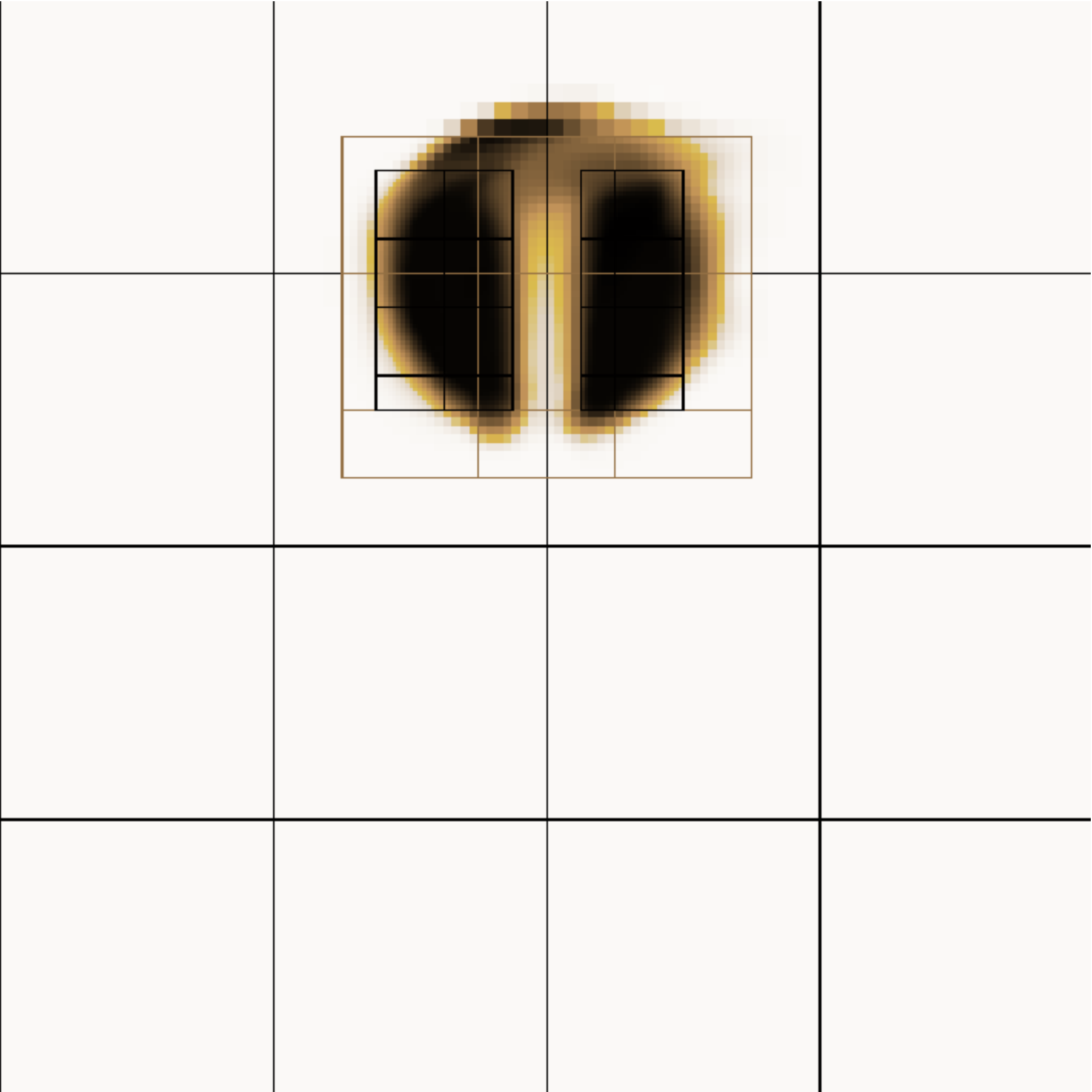}
\caption{GP-AMR prolongation}
\end{subfigure}
\caption{\label{fig:slot1} The morphed slotted cylinder problem at times $t=0.28, 1.44, 2$, from left to right in time.
 Top: (a) Default AMReX with linear prolongation. Bottom: (b) AMReX with adaptive multi-modeled GP prolongation.}
\end{center}
\end{figure}

In this case we wish to find the simulation that best retains 
the profile of the initial condition when it is completed at $t = 2$ as in the previous test. 
We have two levels of refinement on a base grid of size $64 \times 64$ resolution. 
Figure~\ref{fig:slot1} contains snapshots of the simulations at times $t = 0.28, 1.44$ and $t = 2$.  
The goal is to retain as much of the initial condition as possible, in a similar fashion to the 
previous 2D vortex advection test. 

The result shows that the multi-substencil 
GP-AMR prolongation preserves 
the initial condition better than the 
conservative linear scheme native to AMReX.
Notably, there is far less smearing and the circular nature of the cylinder is 
better retained with GP-AMR. 
Furthermore, it should be noted that a larger area of the slotted cylinder is covered by
the finest grid structure
with the GP-AMR prolongation. The refinement criteria in this test and the previous are set for 
critical values of the profile.  This is analogous to refining on regions of high density or pressure. 
We wish to trace the slotted cylinder's evolution with 
the finest grid. 
In this way, we can directly compare the `diffusivity' of each method in how it retains this grid. 
With this test, we see that the default linear prolongation is much more numerically diffusive and 
smears the profile almost immediately, being unable to reconstruct the underlying profile on the coarser grid to the same fidelity that GP-AMR is able to achieve. This results in a far more blurred cylinder at $t=2$. 
While there is some loss with GP-AMR, the profile at $t=2$ far better resembles the cylinder at the 
onset of the simulation.

\subsubsection{Sedov Blast Wave using Castro}
A perhaps more useful test of the algorithm is in a compressible setting, where shock-handling becomes 
necessary. To illustrate the compressible performance of this method, we utilize the Sedov Blast wave
~\cite{sedov}, a radially expanding pressure wave. This simulation is 
solved using Castro with the choice of the piecewise parabolic
method (PPM)~\cite{ppm} for reconstruction along with the Colella and Glaz Riemann 
solver~\cite{ColGlazFerg}. 
For the 2D test we have a base grid of $64 \times 64$, two additional AMR levels, 
using a $r_x=r_y=2$.

\begin{figure}
\begin{center}
\begin{subfigure}{0.49\textwidth}
\includegraphics[width=\textwidth]{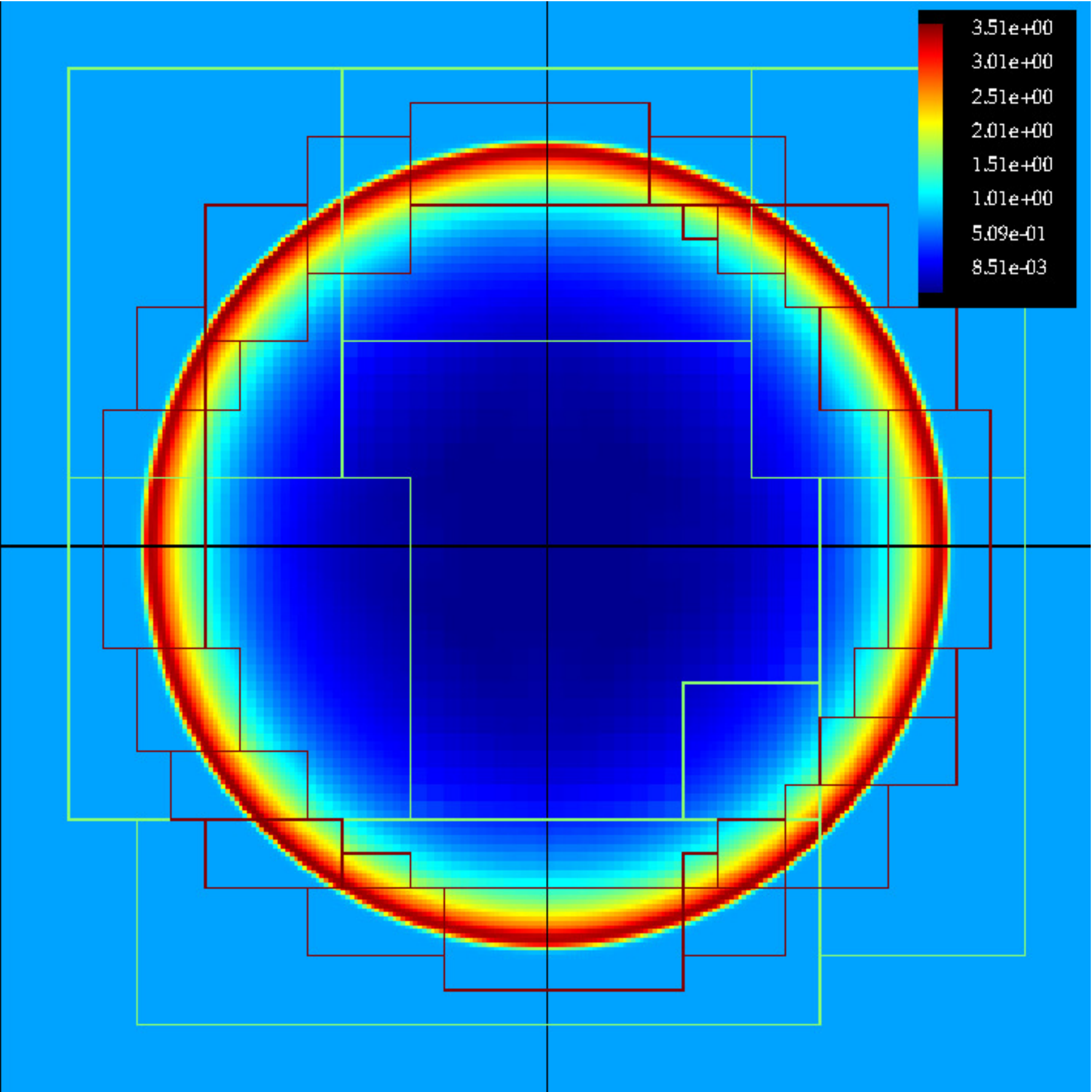}
\caption{Sedov with GP-AMR.}
\end{subfigure}%
\hfill
\begin{subfigure}{0.49\textwidth}
\includegraphics[width=\textwidth]{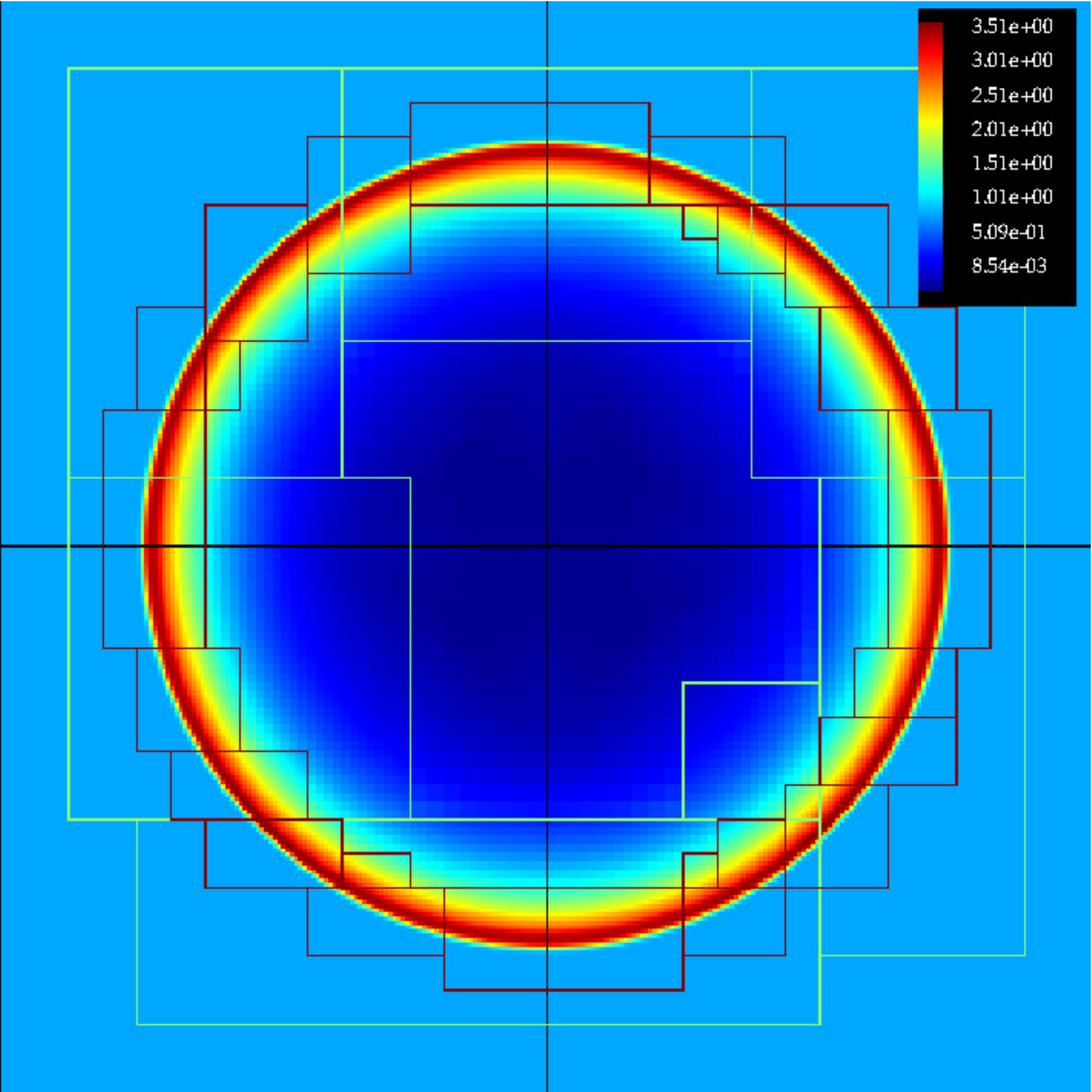}
\caption{Sedov with linear.}
\end{subfigure}
\caption{\label{fig:sedov} A Sedov Blast Wave solution at $t = 0.1$ with two levels of refinement.} 
\end{center}
\end{figure}

Figure~\ref{fig:sedov} illustrates the propagation of the 
Sedov blast wave at $t = 0.1$, and allows us to compare with the linear prolongation method and GP-AMR. 
Although simple, the Sedov blast wave is a good test illustrating the shock-handling capabilities
of the GP multi-substencil model.  Notice that visually, the radial shockwaves in both simulations are 
identical. However, the vacuum in the center of the blast wave is closer to 0 with GP-AMR, as in the 
self-similar analytic solution~\cite{sedov}. 
In Figure~\ref{fig:sedov} 
the AMR levels track the shock as it 
propagates radially and the shock front is contained at the most refined level. At the most refined level,
the shock is handled by the multi-modeled GP-WENO treatment. This increases the computational complexity 
in this region. 
However, the GP algorithm is less expensive in this example, because the shock is very well localized, 
and the majority of the domain is handled by the regular GP model. The standard GP model 
is a simple dot-product using the precomputed weights and the stencil.
 Table~\ref{tab:Sedov} contains the 
performance statistics of the GP-AMR algorithm compared to the default linear using the same workstation
as the previous test. 
 \begin{table}
  \begin{center}
    \caption{Performance of GP-AMR against the default linear AMR on the Sedov Blast Wave	
	      with 6 MPI ranks on an Intel i7-8700K processor.}
    \label{tab:Sedov}
    \begin{tabular}{l|c|c|c} 
       & Execution Time & Prolongation Time & \# of calls\\
	\hline
      $2D$ GP-AMR & 5.719s & 0.07691s & 19743  \\
      $2D$ Linear & 5.698s & 0.12610s & 19439  \\
      $3D$ GP-AMR & 64.27s & 1.28912s & 35202 \\
      $3D$ Non-Adapitve GP-AMR & 72.81s & 5.09467s & 35202 \\ 
      $3D$ Linear & 67.76s & 1.96204s & 35202 \\
     \end{tabular}
  \end{center}
\end{table}

A 3D Sedov blast was also tested, giving us a better look at the multi-substencil cost in the 
shock regions. For this benchmark, the simulation utilized a base grid of 
$32\times 32\times 32$ with additional two levels of AMR, utilizing a refinement factor of 2 
for both levels. 
The wave was advected until $t = 0.01$ with both simulations (GP-AMR and default)
.
Table~\ref{tab:Sedov} also contains the performance metrics for 
the 3D test.  

By setting $\alpha_c=0$ we effectively use the multi-substencil GP-WENO method over all cells. The
metrics for this example are labeled as ``Non-Adaptive GP-AMR'' in Table~\ref{tab:Sedov}. 
Using the multi-substencil GP-WENO method for every grid
 is roughly 5$\times$ more expensive as a prolongation method. This is expected as the multi-modeled 
GP-WENO method combines 5 GP models. 

\subsubsection{Double Mach Reflection using Castro}
The double Mach reflection~\cite{dmr} consists of a Mach 10 shock incident on a reflecting wedge resulting in a complex set of interacting features. 
For this test, we utilize Castro with PPM~\cite{ppm} and the HLLC~\cite{toro} Riemann solver. 
The initial condition describes 
a planar shock front with an angle of $\theta = \pi/3$ extending from the $x$-axis which itself is
a reflecting wall,
\begin{equation}
\begin{pmatrix}
\rho, u, v, p
\end{pmatrix}
=
\begin{cases}
(1.4, \; 0, \; 0, \; 1) & \textrm{for} \quad x > x_{shock}, \\
(8,\; 8.25, \;-8.25, \;116.5) & \textrm{else},
\end{cases}
\end{equation}  
where 
\[x_{shock} = \frac{y + \frac{1}{6}}{\tan{\frac{\pi}{3}}}\]
when $y\in[0,1]$. The full domain of the problem is $[0,4]\times[0,1]$.
 
\begin{figure}
\begin{center}
\includegraphics[width=\textwidth]{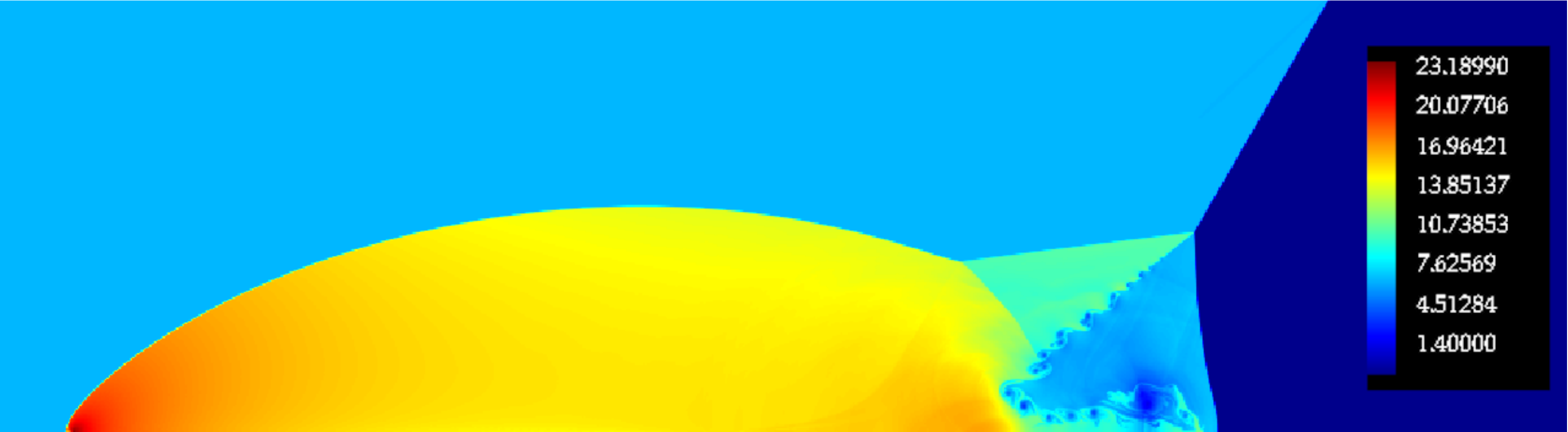}
\caption{\label{fig:dmrfull} The double Mach reflection simulation at $t = 0.2$ with 4 levels of AMR refinement.}
\end{center}
\end{figure}

Figure~\ref{fig:dmrfull} is of the solution to this problem
with 4 levels of refinement starting at a base level with resolution 512$\times$128 
using the GP-based prolongation.

With sufficient resolution and accuracy, we 
observe vortices along the primary slip line, as seen with the copious amount
of vortices in Figure~\ref{fig:dmrfull}. 
The number of vortices serves as a general indication of the numerical diffusivity of the method,
and a quality of Riemann solver. 
In this context, we are interested in the amount of numerical dissipation
of the two different AMR prolongation methods. 

\begin{figure}
\begin{center}
\includegraphics[width=\textwidth]{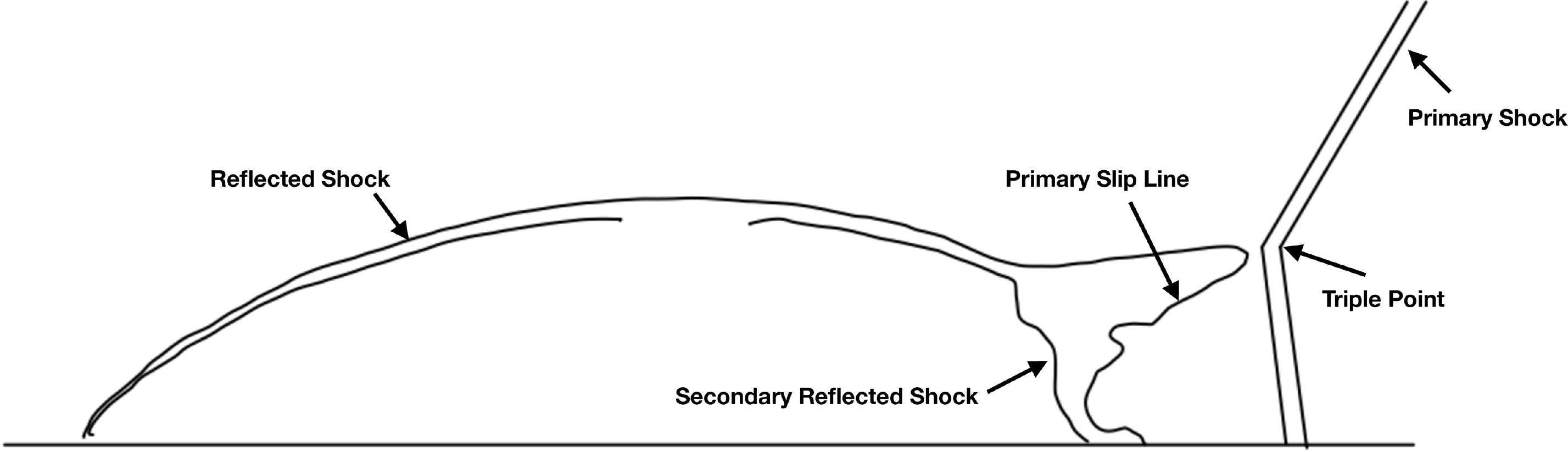}
\caption{\label{fig:dmr_schm} A schematic of the main features in the double mach reflection problem.}
\end{center}
\end{figure}

For reference we present a labeled schematic of the double-mach reflection containing the regions of 
interest for this comparison. Figure~\ref{fig:dmr_schm}
We will refer to features contained in this diagram in the following analysis. Mostly in the central 
region encompassing secondary reflected shock, triple point, and primary slip line. 

\begin{figure}
\begin{center}
\begin{subfigure}{\textwidth}
\centering
\epsfig{width=0.9\textwidth, file=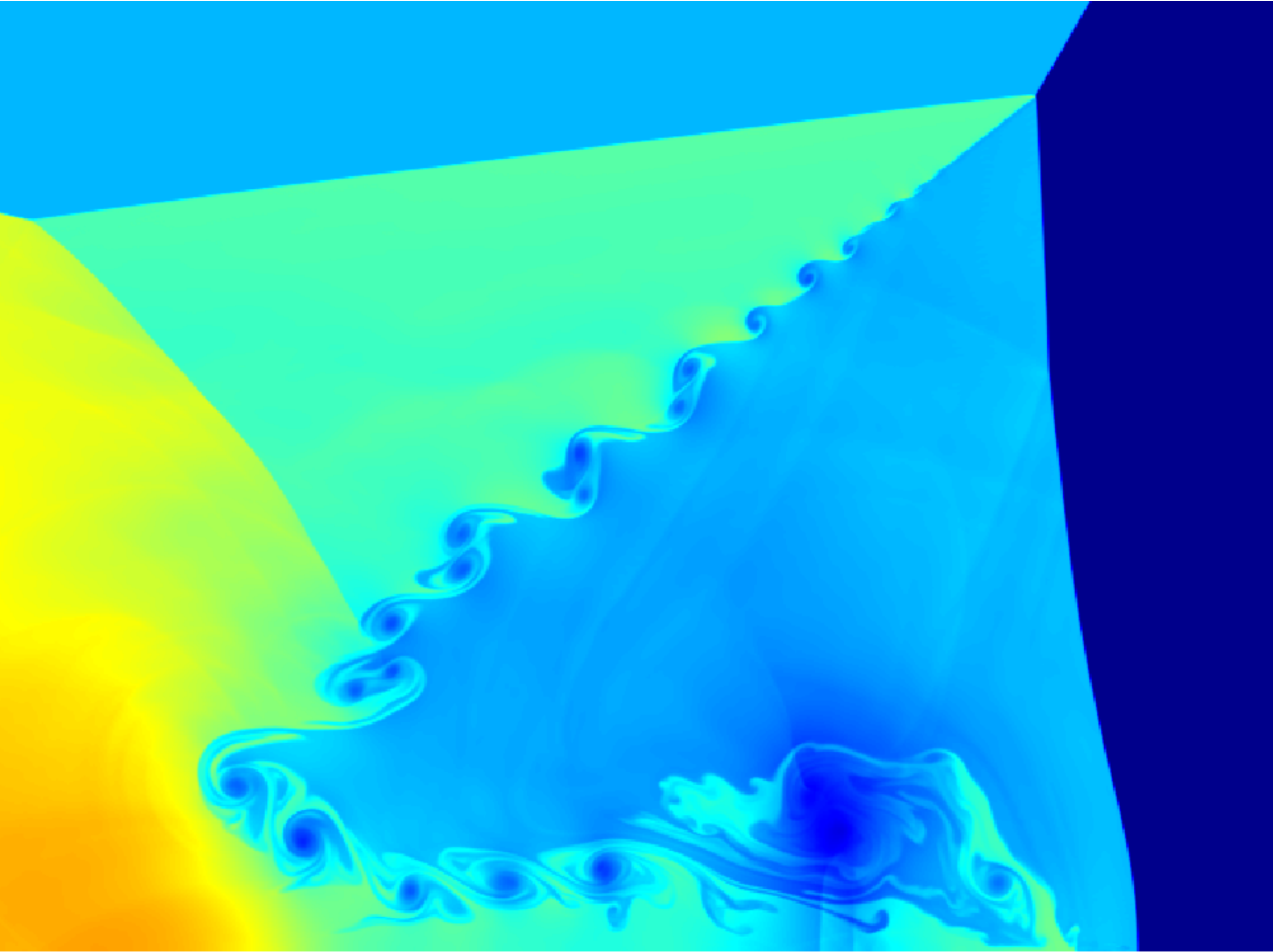}
\caption{\label{fig:GPDMR}}
\bigskip\bigskip
\end{subfigure}
\begin{subfigure}{\textwidth}
\centering
\epsfig{width=0.9\textwidth, file=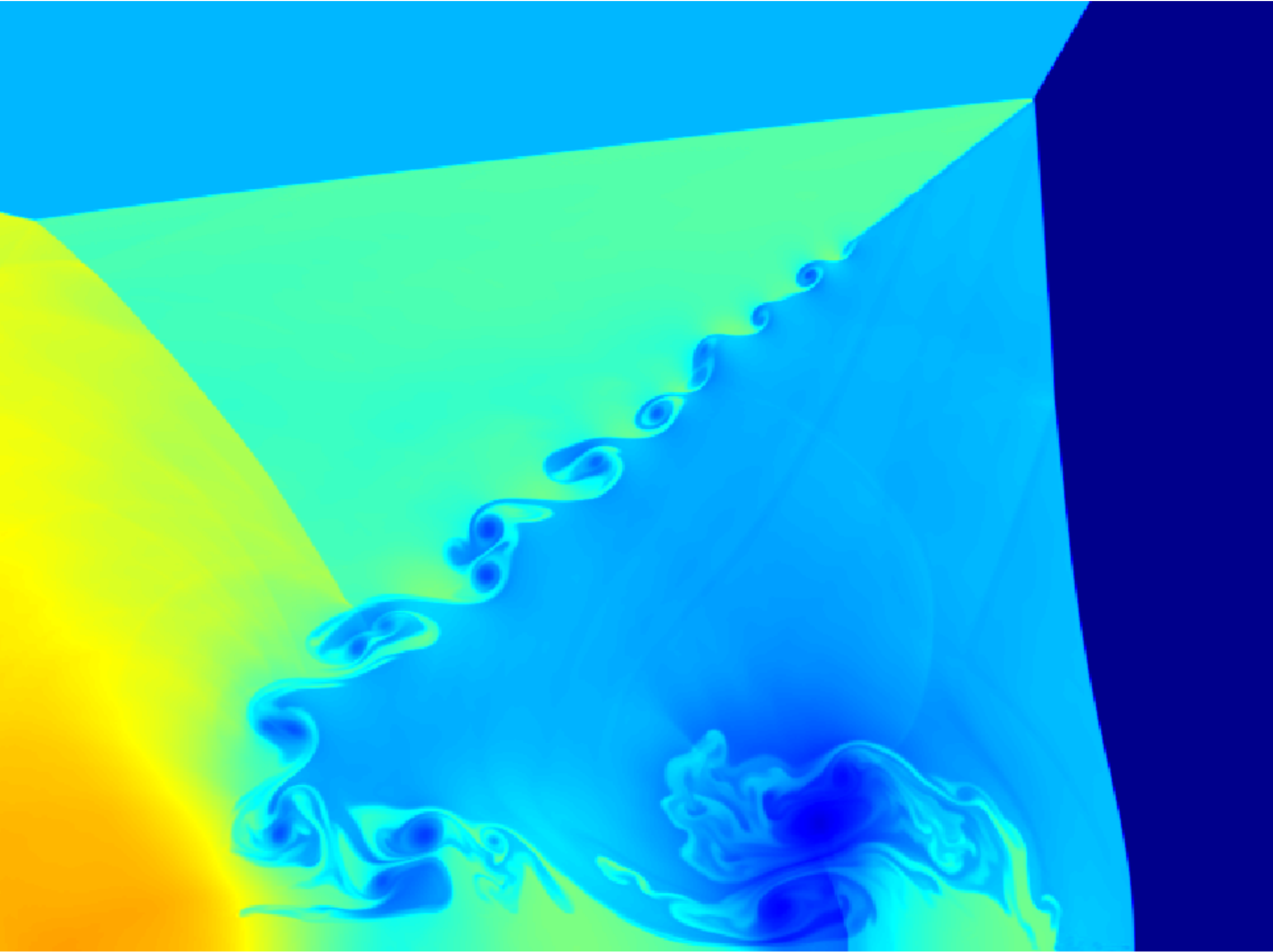}
\caption{\label{fig:LINDMR}}
\end{subfigure}
\caption{\label{fig:dmrComp} (a) GP-AMR simulation, (b) Default linear prolongation based simulation. 
Both simulations visualized in the triple-point region of the domain at $t=0.2$.}
\end{center}
\end{figure}

The default and GP-AMR implementations are compared by zooming into the aforementioned region.  
In Figure~\ref{fig:dmrComp} we observe the effects of the each prolongation method on the
number of vortices along the primary slip line.

As can be seen in Figures~\ref{fig:GPDMR} and~\ref{fig:LINDMR}, 
there is more onset to Kelvin-Helmholtz instability along the primary slip line 
in the GP-based AMR simulation, resulting in more additional vortices above 
the secondary reflected shock wave, 
along with onset to instability on the primary slip line
close to the primary triple point.
As a rudimentary measure, 
the GP-AMR simulation contains 20 vortices in this region whereas
the default AMR contains 17.

With this simulation we see that the default linear prolongation is faster, 
as the adaptive GP algorithm has more cells to prolongate in
the high $\alpha$ regime. Table~\ref{tab:DMR} contains the details about the execution times. 
 \begin{table}
  \begin{center}
    \caption{Performance insights for the Double Mach Reflection problem utilizing 8 nodes and 320 cores
	     on Lux.}
    \label{tab:DMR}
    \begin{tabular}{l|c|c|c} 
       & Execution Time & Prolongation Time & \# of calls\\
	\hline
       GP-AMR & 760.43s & 1.1121s & 39724  \\
       Linear & 705.10s & 0.7395s & 39837  \\
     \end{tabular}
  \end{center}
\end{table}

These results were generated using the University of California, 
Santa Cruz's Lux supercomputer, utilizing 8 nodes. Each node contains two
20-core Intel Xeon Gold 6248 (Cascade Lake) CPUs. This simulation generates more high $\alpha$ regions, 
and thus requires the multi-modeled GP-WENO algorithm more often. This results in the GP-AMR simulation 
being slower than the default prolongation method by 60s. In addition to the increase of computational 
complexity, there is an increase of time spent in the parallel copy algorithm. The multi-modeling 
GP-WENO algorithm requires 2 growth cells at the boarders of each patch, therefore increasing the 
amount of data to be copied.

\subsubsection{Premixed Flame using PeleC}
\begin{figure}
\begin{center}
	\begin{subfigure}{0.49\textwidth}
	\centering{
		\includegraphics[width=.7\textwidth]{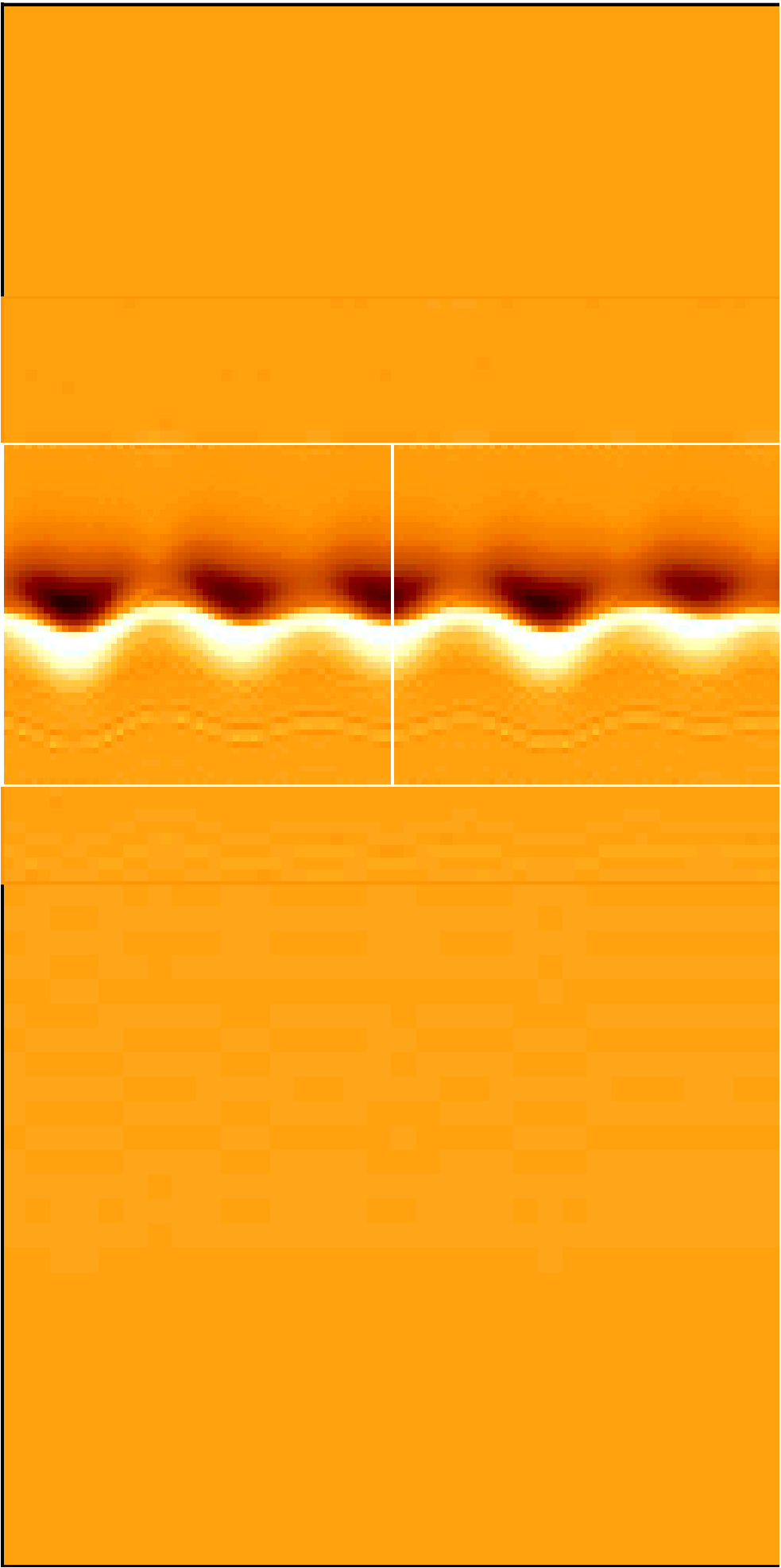}
		\caption{$x-z$ section, $y = 0.3125, z\in[1,4]$}}
	\end{subfigure}
	\bigskip
	\begin{subfigure}{0.49\textwidth}
	\centering{
		\includegraphics[width=.7\textwidth]{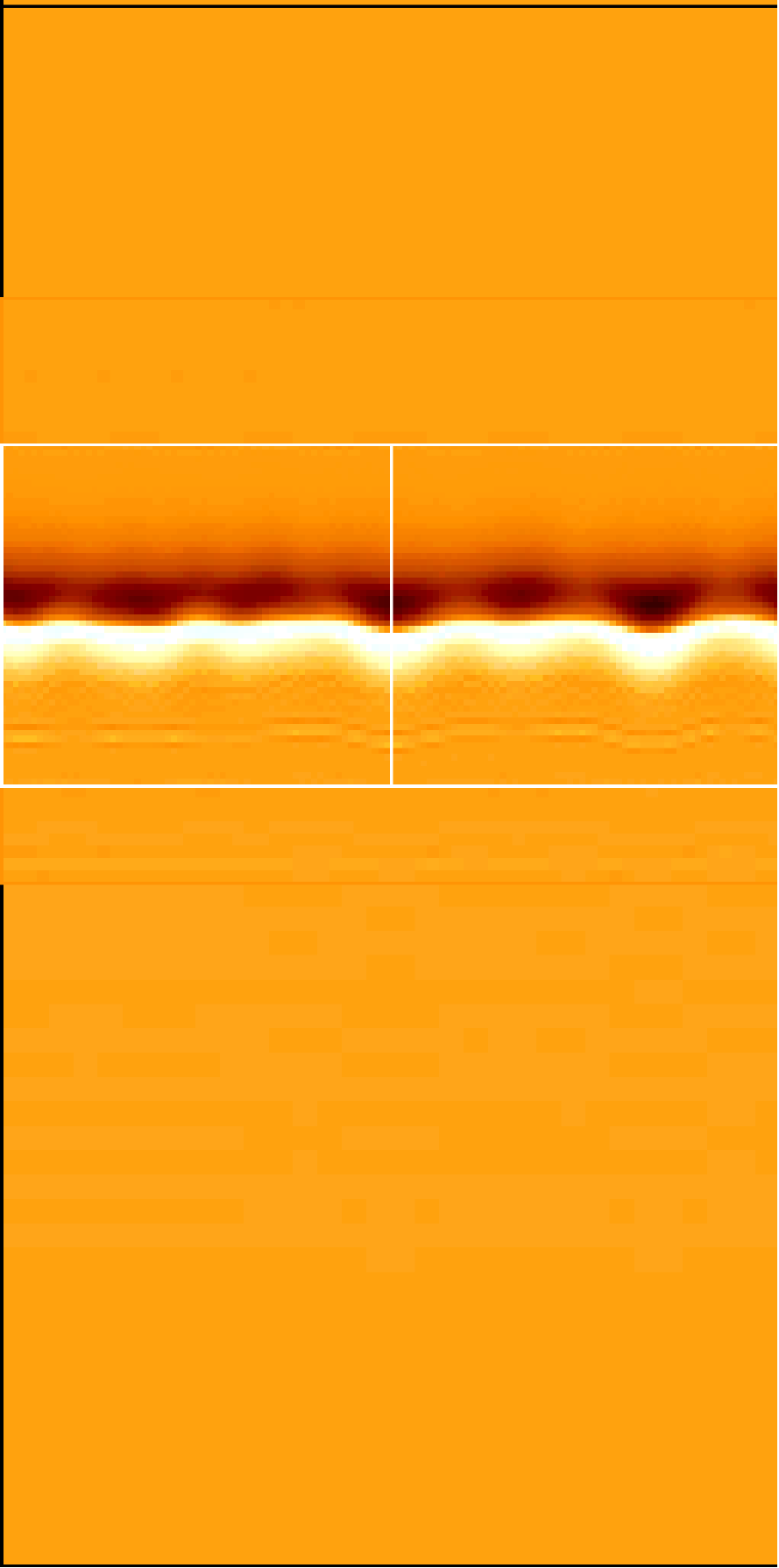}
		\caption{$y-z$ section, $x = 0.3125, z\in[1,4]$}}
	\end{subfigure}
	\begin{subfigure}{0.6\textwidth}
	\centering{
		\includegraphics[width=\textwidth]{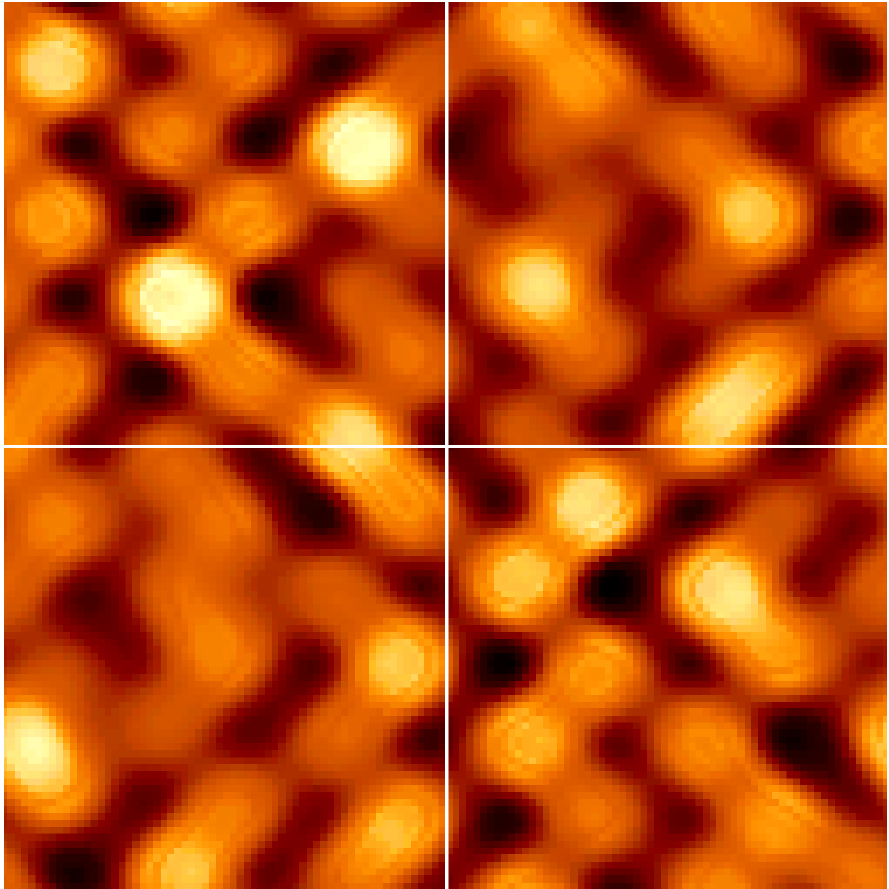}
		\caption{$x-y$ section, $z = 3.0$}}
	\end{subfigure}
	\caption{\label{fig:PMF} Momentum in the $z$-direction of the premixed flame.}
\end{center}
\end{figure}

For the final test problem in this paper,
we produce a steady flame using the AMR compressible combustion simulation code, PeleC~\cite{pelec}. 
With PeleC, chemical species are tracked as mass fractions that are passively moved during the 
advection phase, and diffused subject to transport coefficients and evolved in 
the reaction phase by solving ordinary differential equations  
to compute reaction rates. 
In this problem, the GP-AMR is tied with PPM and the Colella and Glaz Riemann solver~\cite{ColGlazFerg}. We show the 
$\rho w$ (momentum in the $z$ direction) of the premixed flame solution in Figure~\ref{fig:PMF}.
For this illustration the base grid had a $32\times32\times256$ configuration with two additional AMR levels. 
 Furthermore Figure~\ref{fig:PMF} has a color map such that the lighter color are regions with high momentum, 
and the dark regions have low momentum.
 The premixed flame is a 3D flame tube problem 
in a domain that encompasses $[0,0.625]\times[0,0.625]\times[1,6]$. The flame spans the $x$ and $y$ 
dimension and is centered at $z = 3.0$. The gases are premixed and follow Li-Dryer hydrogen
combustion chemical kinetics model~\cite{lidryer}. 

To illustrate some performance metrics we execute the simulation on 32 nodes, 
with 196 NVIDIA V100 GPUs, a base level of 256$\times$128$\times$2048 with two levels of 
refinement. We present the performance of GP-AMR against the default in Table~\ref{tab:perf}. 
We see that the GP-AMR is twice as fast as the default linear on average. 
\begin{table}
  \begin{center}
    \caption{Execution timings of PeleC and AMReX on the Premixed Flame test problem on 32 nodes of the 
Summit supercomputer.}
    \label{tab:perf}
    \begin{tabular}{l|c|c|c} 
       & Execution Time & Prolongation Time & \# of calls \\
      \hline
      GP-AMR & 50.23s & 0.03281s & 940 \\
      Linear & 52.04s & 0.06261s & 940 \\
     \end{tabular}
  \end{center}
\end{table}

\begin{figure}
\begin{center}
\includegraphics[width=0.75\textwidth]{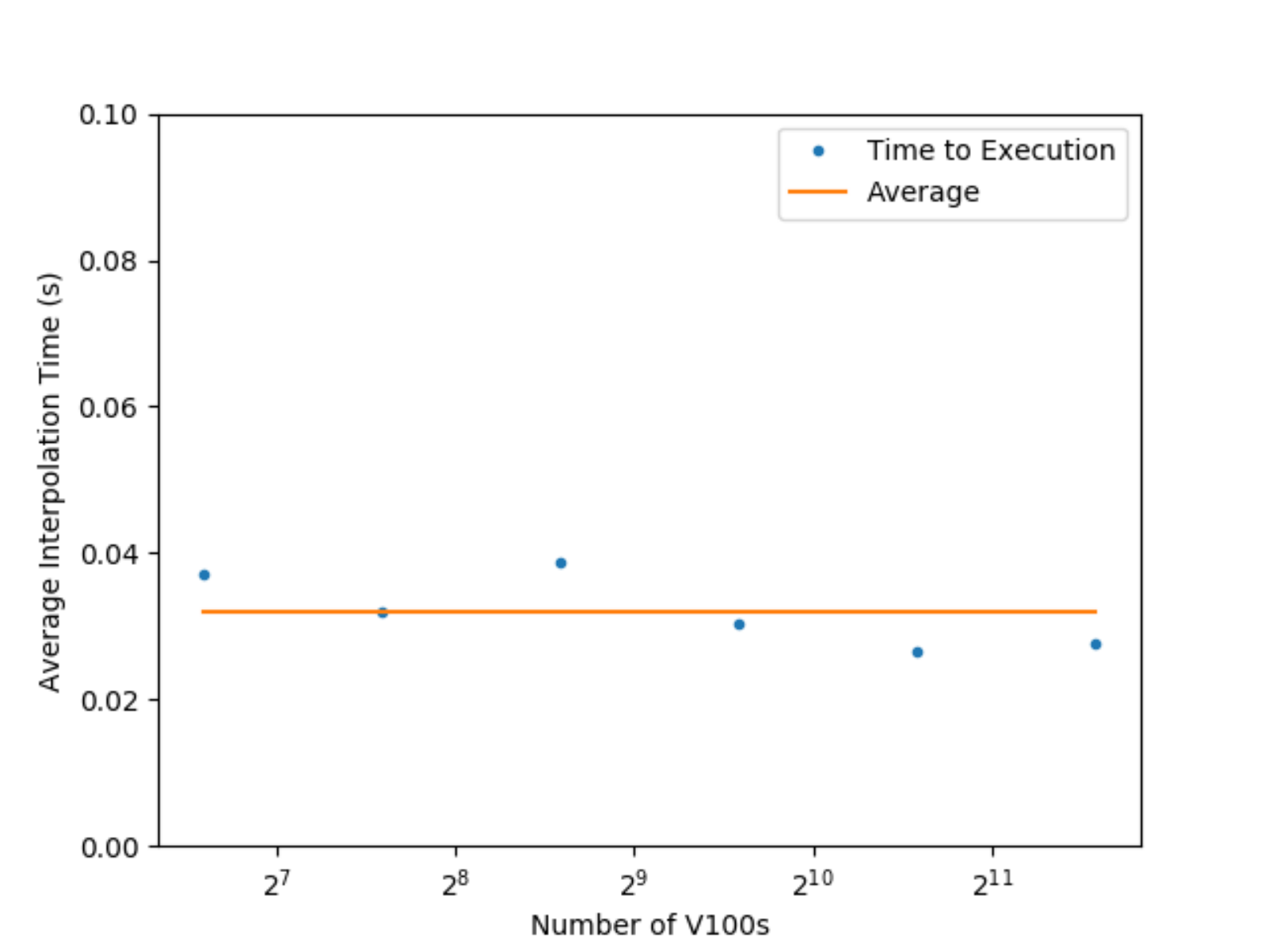}
\caption{\label{fig:wkscl} Weak scaling of GP-AMR utilized in PeleC up to 3072 Nvidia Volta GPUs (512 Nodes on the OLCF Summit Super Computer).} 
\end{center}
\end{figure}

Additionally we perform weak scaling with this problem for up to 3072 NVIDIA V100 GPUs on Summit, with results 
illustrated in Figure~\ref{fig:wkscl}. The $y$-axis of the figure illustrates the average GP prolongation times and the 
$x$-axis is 
number of GPUs from 96 nodes to 3072 GPUs on Summit in logarithmic scale of base-2. Each node on Summit 
contains 6 NVIDIA V100 GPUs, therefore the scaling ranges from 16 to 512 nodes. %
GP-AMR when implemented in AMReX scales very well on Summit, one of the top-class
supercomputers in the modern leadership computing facilities in the world.

\section{Conclusion}
\label{sec:conclusion}
In this paper we developed an efficient, third-order accurate, AMR prolongation method based on
the Gaussian Process Modeling.
This method is general to the type of data being interpolated, 
as illustrated with a substitution of covariance kernels in Eqns.~\ref{eq:sqrexpcov} and~\ref{eq:intsq}. 
In order to handle shock waves, a multi-substencil GP-WENO algorithm inspired by 
WENO~\cite{WENO} was studied. 
We recognize that GP-WENO becomes computationally expensive when shocks are present in simulations.
The tagging approach in Eq.~\eqref{eq:alph} was proposed as a method to 
mitigate this situation. This approach uses 
a grid-scale sized
length-scale parameter, furnished from GP, to detect regions that may contain shocks or non-smooth flows. 
In the three of the five test cases, the 
GP-AMR 
method was faster than the linear interpolation. The other two cases had situations where the patches to be interpolated contained mostly
cells where the 
linear 
GP model did not suffice and required the non-linear multi-substencil treatment. 
Overall, the GP-AMR 
method is a balance between 
speed, stability, and accuracy. In the scope of this paper, the tunable parameters $\ell$ and $\sigma$ are either fixed, or fixed in 
relation to the grid scale. To further adapt the algorithm, one could try and maximize the log of Eqn.~\ref{eq:likely} with 
respect to the hyperparameter $\ell$ as is done in many applications utilizing Gaussian Process regression. However, in our application 
a fixed prescription for $\ell$ appears to hold the properties we desired. The stability of the algorithm is inherently tied to the 
$\sigma$ parameter, which we recommend never being larger than three times the grid scale. In many of the test cases we chose 
$\sigma = 1.5\Delta x$. If additional stability is required, we recommend either tuning $\alpha_c$ to be smaller or to be zero, 
requiring the algorithm to only use the multi-substencil GP model.  

Utilizing the framework provided, an even higher-order prolongation method can be generated by just by increasing the size of the 
stencil while utilizing the same framework. 
This will be inherently useful as more simulations codes are moving to increasingly accurate solutions
with WENO~\cite{WENO,WENO3D} or GP~\cite{reyes_new_2016, reyes2019variable} based reconstruction methods paired with 
Spectral Differed Corrections (SDC)~\cite{sdcode,sdc} which
can yield a fourth or higher order accurate total simulation. In this case, a second order AMR interpolation 
may degrade the overall quality of the solution or incur additional SDC iterations -- increasing the execution time of the simulation.  

\section{Acknowledgements} 
This work was supported in part by the National Science Foundation under grant AST-1908834.
We acknowledge that the current work has come to fruition with help from the Center for Computational Science and Engineering at Lawrence Berkeley National Laboratory, the home of AMReX. 
We thank Dr. Ann S. Almgren, Dr. Weiqun Zhang, and Dr. Marcus Day for their insight and advice when completing this research.  
This research used resources of the Oak Ridge Leadership Computing Facility at the Oak Ridge National Laboratory, which is supported by the Office of Science of the U.S. Department of Energy under Contract No. DE-AC05-00OR22725. 
The first and second authors also acknowledge use of the Lux supercomputer at 
UC Santa Cruz, funded by NSF MRI grant AST 1828315.

\label{sec:references}

\bibliography{mybibfile_merged}

\end{document}